\documentclass{icm2010}
\usepackage{amssymb,latexsym,epsf,graphicx,eepic,bm,color}

\title[Total positivity and cluster algebras]
{Total positivity and cluster algebras}

\author[Sergey Fomin]{Sergey Fomin\thanks{Partially supported by NSF
    grant DMS-0555880.}} 

\contact[fomin@umich.edu]{Department of Mathematics, University of Michigan,
Ann Arbor, MI 48109, USA}

\def\mathcenter#1{\vcenter{\hbox{$#1$}}}
\def\mfigb#1{\mathcenter{\includegraphics[trim=-1 -1 -1 -1]{#1}}}

\newcommand{\darkred}[1]{#1}
\newcommand{\darkblue}[1]{#1}
\newcommand{\darkgreen}[1]{#1}
\newcommand{\lightgreen}[1]{#1}
\newcommand{\lightred}[1]{#1}
\newcommand{\lightblue}[1]{#1}

\numberwithin{equation}{section}

\newtheorem{theorem}{Theorem}[section]

\theoremstyle{definition}

\begin{document}


\begin{abstract}
This is a brief and informal introduction to 
cluster algebras.
It roughly follows
the historical path of their
discovery, made jointly with A.~Zelevinsky. 
Total positivity serves as the main motivation. 
\end{abstract}



\maketitle



\newcommand{\CC}{\mathbb{C}}
\newcommand{\PP}{\mathbb{P}}
\newcommand{\ZZ}{\mathbb{Z}}
\newcommand{\RR}{\mathbb{R}}

\newcommand{\AAA}{\mathcal{A}}
\newcommand{\FFF}{\mathcal{F}}

\newcommand{\LL}{{\mathbf{L}}} 
\newcommand{\xx}{\mathbf{x}}
\newcommand{\zz}{\mathbf{z}}

\newcommand{\Arcs}{\mathbf{A}}
\newcommand{\Exch}{\mathbf{E}}
\newcommand{\Mark}{\mathbf{M}}
\newcommand{\Notch}{\mathbf{N}}
\newcommand{\notch}{\scriptstyle\bowtie} 

\newcommand{\Zmix}{\mathbf{Z}}
\newcommand{\Surf}{\mathbf{S}}

\newcommand{\SM}{(\Surf,\Mark)}

\newcommand{\DDSM}{\Delta^{\hspace{-1pt}\circ}(\Surf,\Mark)}
\newcommand{\DTSM}{\Delta^{\hspace{-1pt}\notch}(\Surf,\Mark)}
\newcommand{\EESM}{\Exch^{\circ\hspace{-1pt}}(\Surf,\Mark)}
\newcommand{\ETSM}{\Exch^{\notch}(\Surf,\Mark)}
\newcommand{\AASM}{\Arcs^{\circ\hspace{-1pt}}(\Surf,\Mark)}
\newcommand{\ATSM}{\Arcs^{\hspace{-1pt}\notch}(\Surf,\Mark)}
\newcommand{\Teich}{\mathcal{T}}
\newcommand{\dTeich}{\widetilde{\Teich}}
\newcommand{\TTSML}{\overline{\Teich}(\Surf,\Mark,\LL)}
\newcommand{\TTSM}{\overline{\Teich}\SM}

\section*{Introduction}

Cluster algebras are encountered in many algebraic and geometric
contexts, with combinatorics providing a unifying framework. 
This short paper reviews the origins of cluster algebras, 
their deep connections with total positivity phenomena, 
and some of their recent manifestations in Teichm\"uller theory. 

The introduction of cluster algebras, made in joint work with
A.~Zelevinsky~\cite{ca1}, 
was rooted in the desire to understand, in a concrete and
combinatorial way, G.~Lusztig's theory of total positivity and canonical
bases in quantum groups (see, e.g., \cite{littelmann, lusztig-survey}). 
Although this goal remains largely elusive (cf.~\cite{leclerc-icm}),
the concept proved valuable due to its surprising ubiquity,
and to the connections it helped uncover between diverse and seemingly
unrelated areas of mathematics. 

This paper gives  a popular and quick introduction 
to the subjects in the title, aimed at an uninitiated reader,
and roughly following the historical order of modern developments
in the two related fields. 
Cumbersome technicalities involved in the usual definition of
cluster algebras are largely omitted, giving way to pro\-to\-typical 
examples from which the reader is invited to generalize, 
to discussions of underlying motivations, and to hints concerning 
further applications and extensions of the basic theory. 
Many important aspects are left out due to space limitations. 

The style is rather informal,
owing to the desire to see the forest through the trees, and to make
the paper accessible to a general mathematical audience.
There are no numbered formulas or theorems:  
results are stated as part of the general narrative.  
Some attributions are missing; they can be found in the sources quoted. 
The goal is to give the reader an intuitive feel for 
what cluster algebras are, and motivate her/him to read the more
  formal expositions elsewhere. 

\pagebreak[3]

Several survey/introductory papers dedicated to the subjects in
the title, approached from various perspectives, have already appeared in the
literature; see~in particular \cite{ando,carter,pcmi,
  tptp, cdm, gls-survey, keller-bourbaki, leclerc-icm, 
zelevinsky-millenium,zelevinsky-sf2001,
  zelevinsky-whatis}. 
An excellent introduction to applications of cluster algebras 
in representation theory is given in B.~Leclerc's
  contribution~\cite{leclerc-icm} to these proceedings. 
Besides consulting these sources and references therein, 
the reader is invited to visit the 
\emph{Cluster Algebras Portal}~\cite{portal}, which 
provides numerous links to publications, conferences, seminars,
thematic programs, 
software packages, etc. 

Our presentation is loosely based on the papers \cite{BFZ, ca3, cats1,
  cats2, fs, tptp, ca1, ca2, cdm},
joint with A.~Berenstein, M.~Shapiro,
D.~Thurston, and A.~Zelevinsky. 
Section~\ref{sec:tp} introduces total positivity and the idea of a
positive/nonnegative part of an algebraic variety. 
Section~\ref{sec:cluster} 
presents the basic notions of cluster algebra theory,
emphasizing its~roots \linebreak[3]
in total positivity. 
Section~\ref{sec:teich} discusses the occurence of cluster algebras in
combinatorial topology of triangulated surfaces, and connections
with Teichm\"uller spaces. 

The format of this brief survey does not allow us to discuss
several important directions of current research on
cluster algebras and related fields. 
In particular, not covered here are the theory of \emph{cluster
  categories} and the various facets of \linebreak[3]
\emph{categorification} \cite{keller-bourbaki,
  keller-categorification, keller-course, nakajima};
the connections between cluster algebras and \linebreak[3]
\emph{Poisson geometry}
\cite{gsv1, gsv2}; 
closely related work 
on cluster
varieties arising in \emph{higher Teichm\"uller theory}
\cite{fg-ihes, fock-gonch}; 
the polyhedral combinatorics of \emph{cluster fans} and \emph{Cambrian
  lattices} \cite{reading-speyer}; 
applications to \emph{discrete integrable systems} 
\cite{k-pdf, yga, iikns, keller-course};
the machinery of \emph{quivers with potentials} \cite{dfz1, dfz2};
connections with \emph{Donaldson-Thomas invariants}
\cite{kontsevich-soibelman, nagao}; 
and other exciting topics.

\smallskip

\noindent
\textbf{Acknowledgments.}
The discovery of cluster algebras, the main work leading~to it,
and the development of fundamentals of the general theory were all
done jointly with my longtime collaborator Andrei Zelevinsky.
I~am indebted to him, and to my co-authors Arkady
Berenstein, Michael Shapiro, and 
Dylan Thurston for their invaluable contributions 
to our joint work discussed below. 
Catharina Stroppel persuaded me to give a talk in Bonn whose design
this presentation follows. 
Bernhard Keller, George Lusztig, and Kelli Talaska made valuable editorial
suggestions.  

\section{Total positivity}
\label{sec:tp}

A matrix $x$ with real entries is called \emph{totally positive} 
(resp., \emph{totally nonnegative}) if all
its minors---that is, determinants of square submatrices---are positive
(resp., nonnegative). 
The first systematic study of these classes of matrices was conducted
in the 1930s by F.~Gantmacher and M.~Krein~\cite{GK},  
following the pioneering work of I.~Schoenberg~\cite{schoenberg}. 
In particular, they 
showed that the eigenvalues of an $n\times n$ totally positive matrix
are real, positive, and distinct. 

Total positivity is a remarkably widespread phenomenon:
matrices with positive/nonnegative minors play an important role
in classical mechanics (theory of small oscillations),
probability (one-dimensional diffusion processes),
discrete potential theory (planar resistor networks),
asymptotic representation theory (the Edrei-Thoma theorem),
algebraic and enumerative combinatorics (immanants, lattice
paths), and of course in linear algebra and its applications. 
See \cite{ando, tptp, GK, gasca-micchelli, karlin} for a plethora of
examples and results, and for additional references. 

\pagebreak[3]

A key technical fact from the classical theory of total positivity 
is C.~Cryer's ``splitting lemma'' \cite{cryer, cryer76}: 
an invertible square matrix $x$ (say of determinant~1) is totally
nonnegative if and only if it has a \emph{Gaussian decomposition}
\[
x=\left[\!\!\begin{array}{ccccc}
1      & 0      & 0      & \cdots & 0      \\
*      & 1      & 0      & \cdots & 0      \\
*      & *      & 1      & \cdots & 0      \\
\vdots & \vdots & \vdots & \ddots & \vdots \\
*      & *      & *      & \cdots & 1 \\
\end{array}\!\right]
\left[\!\!\begin{array}{ccccc}
*      & 0      & 0      & \cdots & 0      \\
0      & *      & 0      & \cdots & 0      \\
0      & 0      & *      & \cdots & 0      \\
\vdots & \vdots & \vdots & \ddots & \vdots \\
0      & 0      & 0      & \cdots & * \\
\end{array}\!\right]
\left[\!\!\begin{array}{ccccc}
1      & *      & *      & \cdots & *      \\
0      & 1      & *      & \cdots & *      \\
0      & 0      & 1      & \cdots & *      \\
\vdots & \vdots & \vdots & \ddots & \vdots \\
0      & 0      & 0      & \cdots & 1 \\
\end{array}\!\right]
\]
in which all three factors 
(lower-triangular unipotent, diagonal, and upper-tri\-an\-gu\-lar unipotent)
are totally nonnegative. 
There is also a counter\-part of this statement for totally positive
matrices. 

The Binet-Cauchy theorem implies that totally positive 
(resp., nonnegative) matrices in $G=\operatorname{SL}_n$ form a
  multiplicative semigroup, denoted by $G_{\ge 0}$. 
In view of Cryer's lemma, the study of $G_{\ge 0}$ can be reduced
to the investigation of its subsemigroup $N_{\ge 0}\subset G_{\ge 0}$ 
of upper-triangular
unipotent totally nonnegative matrices. 

The celebrated Loewner-Whitney Theorem \cite{loewner, whitney}
  identifies the infinitesimal generators of $N_{\ge 0}$ 
as the \emph{Chevalley generators} of the corresponding Lie algebra. 
In pedestrian terms, each upper-triangular unipotent 
totally nonnegative $n\times n$ matrix can be written as a product of
(totally nonnegative) matrices of the form
\[
x_i(t)=
\left[\begin{array}{cccccc}
1      & \cdots & 0      & 0      & \cdots & 0      \\
\vdots & \ddots & \vdots & \vdots & \ddots & \vdots \\
0      & \cdots & 1      & t      & \cdots & 0      \\
0      & \cdots & 0      & 1      & \cdots & 0      \\
\vdots & \ddots & \vdots & \vdots & \ddots & \vdots \\
0      & \cdots & 0      & 0      & \cdots & 1 \\
\end{array}\right];
\]
here the matrix $x_i(t)$ differs from the identity matrix by a single
entry $t\ge0$ in row $i$ and column~$i+1$.
This led G.~Lusztig \cite{lusztig} to the idea of extending the notion
of total positivity to 
other semisimple groups~$G$, by defining the set $G_{\ge 0}$ of
totally nonnegative elements in~$G$ as the semigroup generated by the
Chevalley generators. 
Lusztig has shown that $G_{\ge 0}$ can be 
described by inequalities of the form $\Delta(x)\ge 0$ where
$\Delta$ ranges over the appropriate \emph{dual canonical basis} (at $q=1$).
This set is infinite, and very hard to understand;
fortunately, it can be replaced~\cite{fz-osc} 
by a much simpler (and finite) set of \emph{generalized
  minors}~\cite{fz-dbc}. 

A yet more general (if informal) concept is one of a 
\emph{totally positive/nonnegative variety}. 
Vaguely, the idea is this: take a complex variety~$X$ together with a family
$\mathbf\Delta$ of
``important'' regular functions on~$X$.
The corresponding totally positive (resp., totally nonnegative) 
variety $X_{>0}$ (resp., $X_{\ge 0}$)
is the set of points 
at which all of these functions
take positive (resp., nonnegative) values: 
\[
X_{>0}=\{x\in X:  \Delta(x)>0\ \text{for all $\Delta\in\mathbf\Delta$}\}. 
\]
If $X$ is the affine space of matrices of a given size (or 
$\operatorname{GL}_n(\CC)$ or $\operatorname{SL}_n(\CC)$), 
and $\mathbf\Delta$ is the set of all
minors, then we recover the classical notion. 
One can restrict this construction to matrices lying in a given
stratum of a Bruhat decomposition, or in a given \emph{double Bruhat
  cell}~\cite{fz-dbc, lusztig}. 
Another important example is the \emph{totally positive (resp.,
  nonnegative) Grassmannian}
consisting of the points in a usual Grassmann manifold where all 
Pl\"ucker coordinates can be chosen to be positive (resp.,
nonnegative). 

In each of these examples, the notion of positivity depends on 
a particular choice of a coordinate system: a basis in a vector
space allows us to view linear transformations as matrices; 
a choice of reference flag determines a system of Pl\"ucker
coordinates; and so on. 

Why study totally nonnegative varieties? 
Besides the connections to Lie theory alluded to above, 
there are at least three more reasons. 

First, some totally nonnegative varieties are interesting in their own
right as they can be identified with important spaces, e.g.\ some of 
those arising in \linebreak[3]
Teichm\"uller theory; cf.\ Section~\ref{sec:teich}. 
One can hope to gain additional insight into the structure of such
spaces and their compactifications by
``upgrading'' them to complex varieties, studying associated
quantizations, etc. The nascent ``higher \linebreak[3]
Teichm\"uller theory'' 
\cite{chekhov, fock-gonch} is one prominent expression of this
paradigm. 

Second, passing from a complex variety to its positive part can be
viewed as a step towards its \emph{tropicalization}. 
The deep connections between total positivity, tropical geometry, and
cluster theory lie outside the scope of this short paper;
see \cite{fock-gonch, cats2, ca4} for some aspects of this emerging
research area.

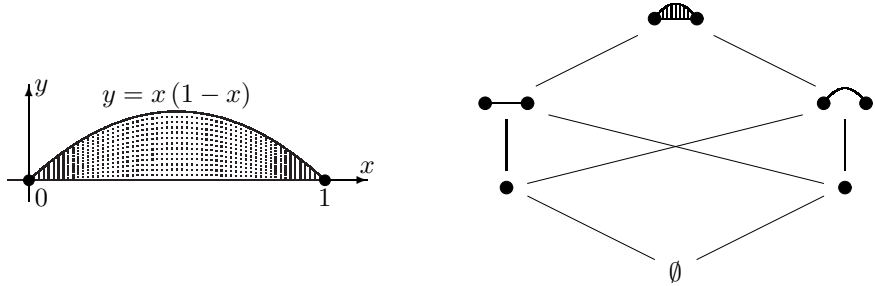
\begin{figure}[ht]
\begin{center}
         \begin{picture}(140,50)(50,-30)
\setlength{\unitlength}{0.8pt} 
\thinlines
\qbezier(0,0)(70,65)(140,0)
\qbezier[60](0,0)(70,60)(140,0)
\qbezier[60](0,0)(70,55)(140,0)
\qbezier[60](0,0)(70,50)(140,0)
\qbezier[60](0,0)(70,45)(140,0)
\qbezier[60](0,0)(70,40)(140,0)
\qbezier[60](0,0)(70,35)(140,0)
\qbezier[60](0,0)(70,30)(140,0)
\qbezier[60](0,0)(70,25)(140,0)
\qbezier[60](0,0)(70,20)(140,0)
\qbezier[60](0,0)(70,15)(140,0)
\qbezier[60](0,0)(70,10)(140,0)
\qbezier[60](0,0)(70,5)(140,0)
\put(-10,0){\line(1,0){150}}
\put(0,0){\circle*{5}}
\put(140,0){\circle*{5}}
\put(70,40){\makebox(0,0){$y=x\,(1-x)$}}
\put(140,0){\vector(1,0){20}}
\put(0,-10){\vector(0,1){55}}
\put(6,-8){\makebox(0,0){$0$}}
\put(140,-8){\makebox(0,0){$1$}}
\put(160,6){\makebox(0,0){$x$}}
\put(6,44){\makebox(0,0){$y$}}
         \end{picture}
\setlength{\unitlength}{1.6pt} 
         \begin{picture}(40,60)(-30,3)
\thinlines
\put(0,0){\makebox(0,0){$\emptyset$}}
\put(40,20){\circle*{3}}
\put(-40,20){\circle*{3}}
\put(-45,40){\circle*{3}}
\put(-45,40){\line(1,0){10}}
\put(-35,40){\circle*{3}}
\put(35,40){\circle*{3}}
\qbezier(35,40)(40,47)(45,40)
\put(45,40){\circle*{3}}
\put(-5,60){\circle*{3}}
\qbezier(-5,60)(0,67)(5,60)
\qbezier[10](-5,60)(0,66)(5,60)
\qbezier[10](-5,60)(0,65)(5,60)
\qbezier[10](-5,60)(0,64)(5,60)
\qbezier[10](-5,60)(0,63)(5,60)
\qbezier[10](-5,60)(0,62)(5,60)
\qbezier[10](-5,60)(0,61)(5,60)
\put(-5,60){\line(1,0){10}}
\put(5,60){\circle*{3}}
\put(-5,3){\line(-2,1){30}}
\put(5,3){\line(2,1){30}}
\put(-35,21){\line(4,1){65}}
\put(35,21){\line(-4,1){65}}
\put(-40,24){\line(0,1){12}}
\put(40,24){\line(0,1){12}}
\put(33,44){\line(-2,1){24}}
\put(-33,44){\line(2,1){24}}
         \end{picture}
\end{center}
\caption{(a) The base $M_{\ge0}$ of the cone $N_{\ge0}$. 
(b) The attachment of algebraic strata.}
\label{fig:digon}
\end{figure}

Yet another reason to study totally nonnegative varieties lies in the
fact that their structure as semialgebraic sets
reveals important features of related complex varieties. 
%
%
We illustrate this phenomenon using the example first studied
in~\cite{lusztig} (cf.~also \cite{BFZ, fs}). 
Consider $N\!\subset\!\operatorname{SL}_n(\CC)$, the subgroup of
$n\!\times\! n$ 
unipotent upper-triangular matrices.
The corresponding totally non\-negative variety is
the semigroup $N_{\ge 0}$ of totally nonnegative matrices in~$N$. 
Take $n=3$; then 
\[
N_{\ge 0}=\left\{
\left[\begin{array}{ccc}
1 & x & y \\
0 & 1 & z \\
0 & 0 & 1 
\end{array}\right]:
\begin{array}{l}
x\ge0\\
y\ge0\\
z\ge0
\end{array}
\ \text{and}\ \ xz-y\ge0
\right\}. 
\]
The inequalities defining $N_{\ge 0}$ are homogeneous in the following
sense: replacing $(x,y,z)$ by $(ax,a^2y, az)$, with $a>0$,  
does not change them. 
Consequently, the space $N_{\ge 0}$ is topologically a
cone with the apex $x\!=\!y\!=\!z\!=\!0$ (the identity matrix) over the base
$M_{\ge 0}\subset N_{\ge 0}$ cut out by the plane $x+z=1$. 
Thus 
\[
M_{\ge 0}\cong \{(x,y)\in\RR^2: 0\le x\le 1\ \text{and}\ \ y\le x(1-x)\}
\]
is the subset of the coordinate plane $\RR^2$ bounded by the $x$ axis
and the parabola $y=x(1-x)$, as shown in Figure~\ref{fig:digon}(a).

The semialgebraic set~$M_{\ge0}$ naturally decomposes into 
5~algebraic strata: two of dimension~$0$, two of dimension~$1$, and
one of dimension~$2$. 
Accordingly, the cone~$N_{\ge0}$ decomposes into 6~algebraic
strata of dimension 1 higher; the apex of~$N_{\ge0}$ corresponds to
the ``empty face'' of~$M_{\ge0}$.  
See 
Figure~\ref{fig:digon}(b). 

\pagebreak[3]

The adjacency of these strata is described by a partial order
isomorphic~to the 
\emph{Bruhat order} on the symmetric group~$\mathcal{S}_3\,$. 
This happens in general, for any~$n$:~the decomposition of $N_{\ge0}$
into algebraic strata produces a 
CW-complex with cell attachments described by the Bruhat order
on~$\mathcal{S}_n\,$. 
Recall that the same partial order describes the attachment of Schubert
cells in the manifold of complete flags in~$\CC^n$. 
The latter has rich topology, and is a central object of study in
modern Schub\-ert Calculus. 
By contrast, $N_{\ge0}$ and $M_{\ge0}$ have no topology to speak of 
(in fact, $M_{\ge0}$ is expected to be homeomorphic to a
ball~\cite{fs, hersh})
but has a cell decomposition with exactly the same cell attachments. 
The big difference of course is that the complex
Schubert cells have twice the dimensions of their real (more
precisely, positive real) counterparts. 
Still, the stratification of $M_{\ge0}$ resulting from its
semialgebraic structure somehow ``remembers''the Bruhat order---which is all
one needs to know in order to reconstruct the topology of 
the flag manifold and its Schubert cells/varieties---including Schubert
and Kazhdan-Lusztig polynomials,~etc. 

\section{Cluster algebras}
\label{sec:cluster}

The discussion in Section~\ref{sec:tp} 
prompts one to ask:
Which algebraic varieties $X$ 
have a natural notion of positivity? 
Which families $\mathbf\Delta$ of regular functions should one consider
in defining this notion? 
The concept of a cluster algebra can be viewed as an attempt to 
provide a general answer to these questions. 
Since the definition 
is fairly technical, 
we start with an example and then generalize. 

Our prototypical example of a cluster algebra $\AAA$ is 
the coordinate ring of the
\emph{base affine space} for the special linear group 
$G=\operatorname{SL}_n(\CC)$, defined as follows. \linebreak[3]
The subgroup  $N\subset G$ 
of unipotent upper-triangular matrices acts on
$G$ by right multiplication. 
The algebra $\AAA=\CC[G/N]$ 
consists of regular functions on~$G$ 
which are invariant under this action~of~$N$. 
Thus elements of $\AAA$ can be viewed as polynomials in the entries $x_{ij}$
of a matrix $x\!=\!(x_{ij})\!\in\!\operatorname{SL}_n(\CC)$ 
which are invariant under column
operations that add to a column of~$x$ a linear combination of
preceding columns. 
%
%
Classical invariant theory tells us that $\AAA$ is
generated by the 
\emph{flag minors} 
\[
\Delta_I: x 
\mapsto \det(x_{ij}|i\in I, j\le|I|)   
\]
where $I$ ranges over nonempty proper subsets of $\{1,\dots,n\}$. 
That is, $\Delta_I$ is a minor occupying the rows in~$I$ and the first
several columns. 
The generators $\Delta_I$ satisfy certain well known homogeneous quadratic
identities sometimes called 
\emph{generalized Pl\"ucker relations}. 

A point in $G/N$ represented by a matrix $x$ is, by
definition, totally positive/nonnegative if all flag minors
$\Delta_I$ take positive/nonnegative values~at~$x$. 
Total positivity in $G/N$ is closely related to the classical
notion of total positivity in~$G$:
it is not hard to deduce from Cryer's lemma
that a matrix $x$ is totally positive if and only if both $x$ and its
transpose represent totally positive elements in~$G/N$. 

There are $2^n-2$ flag minors;
do we really have to test all of them to verify that a point $x\in
G/N$ is totally positive?
The answer is no: it suffices to test positivity of
$\dim(G/N)=\frac{(n-1)(n+2)}{2}$ minors; 
one could hardly hope for a more efficient~test. 

\begin{figure}[ht]
\begin{center}
\setlength{\unitlength}{2.4pt} 
\begin{picture}(60,60)(0,3) 
\thinlines 

\put(0,0){\line(1,1){60}}
\put(0,20){\line(1,-1){15}}
\put(0,40){\line(1,-1){35}}
\put(60,40){\line(-1,-1){35}}
\put(60,20){\line(-1,-1){15}}
\put(0,60){\line(1,-1){60}}

\qbezier(35,5)(40,0)(45,5)
\qbezier(15,5)(20,0)(25,5)

\put(-3,0){\makebox(0,0){$\mathbf{1}$}}
\put(-3,20){\makebox(0,0){$\mathbf{2}$}}
\put(-3,40){\makebox(0,0){$\mathbf{3}$}}
\put(-3,60){\makebox(0,0){$\mathbf{4}$}}

\put( 3,10){\makebox(0,0){$\Delta_{1}$}}
\put(20,10){\makebox(0,0){$\Delta_{2}$}}
\put(40,10){\makebox(0,0){$\Delta_{3}$}}
\put(57,10){\makebox(0,0){$\Delta_{4}$}}

\put(3,30){\makebox(0,0){$\Delta_{12}$}}
\put(30,20){\makebox(0,0){$\Delta_{23}$}}
\put(57,30){\makebox(0,0){$\Delta_{34}$}}

\put(3,50){\makebox(0,0){$\Delta_{123}$}}
\put(57,50){\makebox(0,0){$\Delta_{234}$}}

\end{picture}
\qquad\qquad
\begin{picture}(60,60)(0,3) 
\thinlines 

\put(0,0){\line(1,1){60}}
\put(0,40){\line(1,-1){35}}
\put(60,40){\line(-1,-1){20}}
\put(60,20){\line(-1,-1){15}}
\put(0,60){\line(1,-1){60}}
\put(20,30){\line(1,-1){10}}

\qbezier(30,20)(35,15)(40,20)

\qbezier(35,5)(40,0)(45,5)

\qbezier(10,30)(15,35)(20,30)
\put(0,20){\line(1,1){10}}

\put(-3,0){\makebox(0,0){$\mathbf{1}$}}
\put(-3,20){\makebox(0,0){$\mathbf{2}$}}
\put(-3,40){\makebox(0,0){$\mathbf{3}$}}
\put(-3,60){\makebox(0,0){$\mathbf{4}$}}

\put( 3,10){\makebox(0,0){$\Delta_{1}$}}
\put(20,26){\makebox(0,0){$\Delta_{13}$}}
\put(40,10){\makebox(0,0){$\Delta_{3}$}}
\put(57,10){\makebox(0,0){$\Delta_{4}$}}

\put(3,30){\makebox(0,0){$\Delta_{12}$}}
\put(31,25){\makebox(0,0){$\Delta_{23}$}}
\put(57,30){\makebox(0,0){$\Delta_{34}$}}

\put(3,50){\makebox(0,0){$\Delta_{123}$}}
\put(57,50){\makebox(0,0){$\Delta_{234}$}}

\end{picture}

\end{center}
\caption{Two pseudoline arrangements, and associated chamber minors}
\label{fig:4-pseudo}
\end{figure}
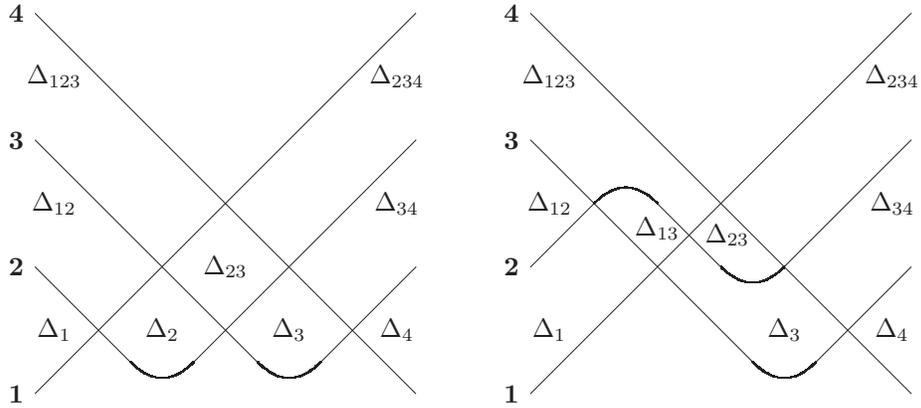

To design such tests, 
we will need the notion of a \emph{pseudoline arrangement}. \linebreak[3]
The latter is a collection of $n$ ``pseudolines'' 
each of which is a graph of a continuous function on $[0,1]$; 
each pair of pseudolines must have exactly one crossing point in
common. 
(See Figure~\ref{fig:4-pseudo}.)
The resulting arrangement is considered up to isotopy. 

We label the pseudolines $\mathbf{1}$
  through~$\mathbf{n}$ by numbering their left endpoints from the bottom up. 
To each \emph{region} $R$ of a pseudoline arrangement, 
with the exception of the top and the bottom regions,
we associate the \emph{chamber minor} $\Delta_{I(R)}$ 
(cf.~\cite{BFZ}) 
defined as the flag minor indexed by the set $I(R)$ of labels
of the pseudolines passing \emph{below}~$R$. 
The $\frac{(n-1)(n+2)}{2}$ chamber minors associated
with a given pseudoline arrangement form an \emph{extended cluster}; 
we shall see that the positivity of these minors implies 
that \emph{all} flag minors of a given matrix are positive. 

There are two types of regions: 
the \emph{bounded} regions entirely surrounded by 
pseudolines, and the \emph{unbounded} ones, 
adjacent to the left and right borders. 
The $2(n\!-\!1)$ chamber minors associated with unbounded regions are called
\emph{frozen}: 
these minors are present in every arrangement. 
For $n=4$, the frozen minors are 
$\Delta_1, \, \Delta_{12}, \, \Delta_{123}, \, \Delta_4, \,
  \Delta_{34}, \,$ and $\Delta_{234}$ 
(cf.\ Figure~\ref{fig:4-pseudo}).

The chamber minors corresponding to the bounded
regions form the \emph{cluster} associated with the
given pseudoline arrangement. 
(Thus an extended cluster is a cluster plus the frozen
minors.) 
Each cluster contains $\binom{n-1}{2}$
chamber minors. 
The two pseudoline arrangements shown in
Figure~\ref{fig:4-pseudo} have clusters
$\{\Delta_2,\Delta_3,\Delta_{23}\}$  and 
$\{\Delta_{13},\Delta_3,\Delta_{23}\}$, respectively.

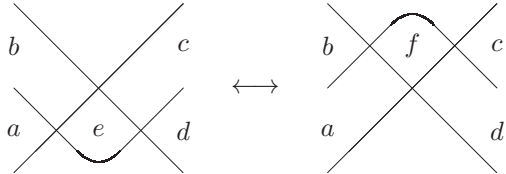
\begin{figure}[ht]
\begin{center}
\setlength{\unitlength}{1.6pt} 
\begin{picture}(40,40)(0,3) 
\thinlines 
\put(0,0){\line(1,1){40}}
\put(0,40){\line(1,-1){40}}
\put(0,20){\line(1,-1){15}}
\put(40,40){\line(-1,-1){35}}
\put(40,20){\line(-1,-1){15}}
\qbezier(15,5)(20,0)(25,5)
\put(0,10){\makebox(0,0){$a$}}
\put(0,30){\makebox(0,0){$b$}}
\put(40,30){\makebox(0,0){$c$}}
\put(40,10){\makebox(0,0){$d$}}
\put(20,10){\makebox(0,0){$e$}}
\end{picture}
\begin{picture}(30,40)(0,3) 
\put(15,20){\makebox(0,0){$\longleftrightarrow$}}
\end{picture}
\begin{picture}(40,40)(0,3) 
\thinlines 
\put(0,0){\line(1,1){40}}
\put(0,40){\line(1,-1){40}}
\put(0,20){\line(1,1){15}}
\put(40,40){\line(-1,-1){35}}
\put(40,20){\line(-1,1){15}}
\qbezier(15,35)(20,40)(25,35)
\put(0,10){\makebox(0,0){$a$}}
\put(0,30){\makebox(0,0){$b$}}
\put(40,30){\makebox(0,0){$c$}}
\put(40,10){\makebox(0,0){$d$}}
\put(20,30){\makebox(0,0){$f$}}
\end{picture}
\end{center}
\caption{A local move in a pseudoline arrangement}
\label{fig:move}
\end{figure}

These two clusters differ in one element only. 
This is because the corresponding two arrangements 
are related to each other by a \emph{local move} 
consisting in dragging one of the pseudolines 
through an intersection of two others; see Figure~\ref{fig:move}. 
As a result of such a move, one chamber minor (namely~$e$ in
Figure~\ref{fig:move}, and $\Delta_2$ in Figure~\ref{fig:4-pseudo}) 
disappears (we say that this minor is \emph{flipped}), and a new
one (namely~$f$ in Figure~\ref{fig:move}, 
and $\Delta_{13}$ in Figure~\ref{fig:4-pseudo}) is introduced.

It can be shown that for a local move as in Figure~\ref{fig:move},
the chamber minors associated with the regions where the action takes
place satisfy the identity
 \[
ef=ac+bd. 
\]
This identity is one of the generalized Pl\"ucker relations
alluded to above. 
We call it an \emph{exchange relation}, 
as the chamber minors~$e$ and~$f$ are exchanged
by the local move. 
For the local move shown in Figure~\ref{fig:4-pseudo}, the exchange
relation~is 
\[
\Delta_2
\Delta_{13}={\Delta_{12}}\,\Delta_3+{\Delta_1}\,\Delta_{23}\,. 
\]

The new chamber minor $f$ produced by a local move is given by a
simple rational expression $f=\frac{ac+bd}{e}$
in the chamber minors of the original arrangement.  
Note that this expression is \emph{subtraction-free} (no minus
signs). 
One can now start with a particular pseudoline arrangement,
label its regions by indeterminates, then use iterated local moves 
(combined with the corresponding birational
transformations) to generate all possible arrangements, and
in doing so write \emph{all} flag minors as rational expressions in
the initial extended cluster.  
All these expressions are clearly subtraction-free, and the claim
follows:
if the elements of the initial extended cluster evaluate positively at
a given point in~$G/N$, then so do all flag minors. 

Let $\FFF$ denote the field of rational functions in the formal
variables making up the initial extended  cluster.
Inside~$\FFF$, 
the rational expressions discussed in the previous paragraph generate 
the subalgebra~$\AAA$ canonically isomorphic to~$\CC[G/N]$. 
Notice that our construction does not explicitly involve the 
group~$G$: we can pretend to be unaware that we are dealing with matrices,
their minors, etc. \linebreak[3]
Yet the construction produces, by design, an algebra~$\AAA$ 
equipped with a distinguished set of generators~$\mathbf{\Delta}$ (the
rational expressions corresponding to the flag minors), 
and thus endowed with a notion of (total) positivity. 

The example of a base affine space treated above 
displays, in a rudimentary form, 
the main features of a general cluster algebra set-up.
We next proceed to describing the latter on an informal level, with
details to be filled in later on. 

Fix a field $\FFF$ of rational functions in several variables, some of
which are designated as ``frozen.'' 
Imagine a (potentially infinite) family of
equinumerous finite collections (``clusters'') of 
elements in~$\FFF$. 
(These elements, called \emph{cluster variables}, 
can be thought of as regular functions on some ``cluster variety''~$X$.) 
Each cluster can be ``extended'' by adjoining the frozen
variables. 
The (extended) clusters are the vertices of a connected regular 
graph 
in which adjacent clusters are related by
birational transformations of the most simple
kind, replacing an arbitrary element of a cluster by a sum of two
monomials divided by the element being removed. \linebreak[3]
(By a monomial we mean a product of elements of a given extended
cluster.) 
These transformations are subtraction-free, so 
positivity of the elements of a cluster at a point $x\in X$
does not depend on the choice of a cluster. 
The birational maps between adjacent clusters are encoded by
appropriate combinatorial data, and 
the construction is made rigid by 
mandating that these data are transformed (as one moves to an adjacent
cluster) according to certain canonical rules. 
These combinatorial rules define a discrete dynamics that drives
the algebraic dynamics of cluster transformations. 
Consequently, the choice of initial combinatorial data
(the pseudoline arrangement in the example of~$G/N$)
determines, in a recursive fashion, the entire structure of clusters
and exchanges. 
The corresponding cluster algebra is then defined as the subring of
the ambient field~$\FFF$ generated by the elements of all extended
clusters. 

In the example of the base affine space,
one key feature of the set-up described above is lacking: 
we do not always know how to exchange an element of a cluster. \linebreak[3]
If a region in a pseudoline arrangement 
is bounded by more than three pseudolines, then the
corresponding chamber minor cannot be readily flipped by a local move. 
For instance, 
how do we exchange the chamber minor~$\Delta_{23}$ in
Figure~\ref{fig:4-pseudo} on the left? 
There is in fact a ``hidden'' exchange relation of
the form 
$\Delta_{23}\,\,\circledast=\circledast+\circledast$
---but how do we guess what those $\circledast$'s are? 

The answer to this question will fall into our lap 
once we replace the language of pseudoline
arrangements, too specialized for a general theory, 
by a more universal combinatorial language of quivers. 
(Using quivers somewhat restricts the generality of the cluster theory,
but is general enough for the purposes of this paper.) 
Developing this language will take a little time---but will 
pay off quickly.  

A \emph{quiver} is a finite oriented graph.
We allow multiple edges, but not loops 
(i.e., edges connecting a vertex to itself) 
or oriented 2-cycles 
(i.e.,~edges of opposite orientation connecting the same pair of
vertices). 
We will need a slightly richer notion, with some
vertices in a quiver designated as \emph{frozen}.
The remaining vertices are called \emph{mutable}. 
We assume that no edges connect frozen vertices to each other. 
(Such edges would make no difference in what follows.) 

Quivers play the role of the aforementioned
combinatorial data accompanying the clusters. 
We think of the vertices of a quiver as labeled by the elements of an
extended cluster, 
so that the frozen vertices are labeled by the
frozen variables, and the mutable
vertices by the cluster variables.  

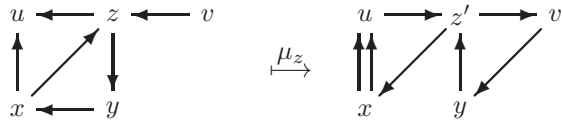
\begin{figure}[ht]
\begin{center}
\setlength{\unitlength}{1.8pt} 
\begin{picture}(40,15)(0,7) 

\put( 0,0){\makebox(0,0){$x$}}
\put(20,0){\makebox(0,0){$y$}}

\put(0,20){\makebox(0,0){$u$}}
\put(20,20){\makebox(0,0){$z$}}
\put(40,20){\makebox(0,0){$v$}}

\thicklines 

\put(16,0){\vector(-1,0){12}}
\put(16,20){\vector(-1,0){12}}
\put(36,20){\vector(-1,0){12}}

\put(20,16){\vector(0,-1){12}}
\put(0,4){\vector(0,1){12}}

\put(3,3){\vector(1,1){14}}

\end{picture}
\qquad$\stackrel{\displaystyle\mu_z}{\longmapsto}$\qquad
\begin{picture}(40,15)(0,7) 

\put( 0,0){\makebox(0,0){$x$}}
\put(20,0){\makebox(0,0){$y$}}

\put(0,20){\makebox(0,0){$u$}}
\put(20,20){\makebox(0,0){$z'$}}
\put(40,20){\makebox(0,0){$v$}}

\thicklines

\put(4,20){\vector(1,0){12}}
\put(24,20){\vector(1,0){12}}

\put(-1.3,4){\vector(0,1){12}}
\put(1.3,4){\vector(0,1){12}}
\put(20,4){\vector(0,1){12}}

\put(17,17){\vector(-1,-1){14}}
\put(37,17){\vector(-1,-1){14}}

\end{picture}
\end{center}
\caption{A quiver mutation. Vertices $u$ and $v$ are frozen.}
\label{fig:quiver-mutation}
\end{figure}

We next describe the quiver analogue 
of a local move. 
Let $z$ be a mutable vertex in a quiver~$Q$.
The \emph{quiver mutation} $\mu_z$ transforms $Q$ into a new
quiver~$Q'=\mu_z(Q)$ via a sequence of three steps.
At the first step, for each pair of directed edges $x\to z\to y$
passing through~$z$, we introduce a new edge $x\to y$ (unless both
$x$ and~$y$ are frozen, in which case do nothing). 
At the second step, we reverse the direction of all edges incident
to~$z$. 
At the third step, we repeatedly remove oriented 2-cycles
until unable to do so. See Figure~\ref{fig:quiver-mutation}. 
It is easy to check that mutating $Q'$ at $z'$ recovers~$Q$.  

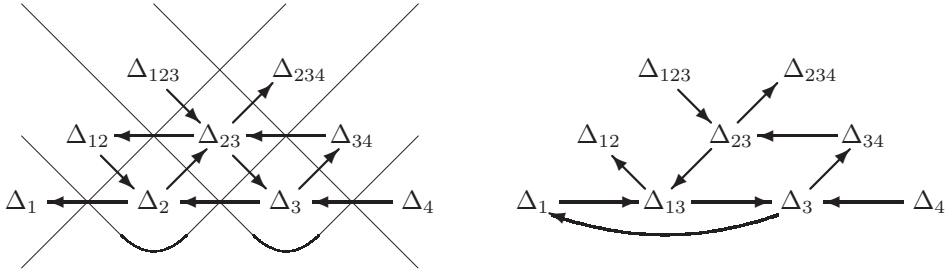
\begin{figure}[ht]
\begin{center}
\setlength{\unitlength}{2.5pt} 
\begin{picture}(60,35)(0,6) 
\thinlines 

\put(0,0){\line(1,1){40}}
\put(0,20){\line(1,-1){15}}
\put(0,40){\line(1,-1){35}}
\put(60,40){\line(-1,-1){35}}
\put(60,20){\line(-1,-1){15}}
\put(20,40){\line(1,-1){40}}

\qbezier(35,5)(40,0)(45,5)
\qbezier(15,5)(20,0)(25,5)

\put( 0,10){\makebox(0,0){$\Delta_{1}$}}
\put(20,10){\makebox(0,0){$\Delta_{2}$}}
\put(40,10){\makebox(0,0){$\Delta_{3}$}}
\put(60,10){\makebox(0,0){$\Delta_{4}$}}

\put(10,20){\makebox(0,0){$\Delta_{12}$}}
\put(30,20){\makebox(0,0){$\Delta_{23}$}}
\put(50,20){\makebox(0,0){$\Delta_{34}$}}

\put(20,30){\makebox(0,0){$\Delta_{123}$}}
\put(42,30){\makebox(0,0){$\Delta_{234}$}}

\thicklines 

\put(16,10){\vector(-1,0){12}}
\put(36,10){\vector(-1,0){12}}
\put(56,10){\vector(-1,0){12}}

\put(26,20){\vector(-1,0){12}}
\put(46,20){\vector(-1,0){12}}

\put(22,12){\vector(1,1){6}}
\put(42,12){\vector(1,1){6}}
\put(32,22){\vector(1,1){6}}

\put(12,17){\vector(1,-1){5}}
\put(32,17){\vector(1,-1){5}}
\put(22,27){\vector(1,-1){5}}

\end{picture}
\qquad\qquad
\begin{picture}(60,35)(0,6) 

\put( 0,10){\makebox(0,0){$\Delta_{1}$}}
\put(20,10){\makebox(0,0){$\Delta_{13}$}}
\put(40,10){\makebox(0,0){$\Delta_{3}$}}
\put(60,10){\makebox(0,0){$\Delta_{4}$}}

\put(10,20){\makebox(0,0){$\Delta_{12}$}}
\put(30,20){\makebox(0,0){$\Delta_{23}$}}
\put(50,20){\makebox(0,0){$\Delta_{34}$}}

\put(20,30){\makebox(0,0){$\Delta_{123}$}}
\put(42,30){\makebox(0,0){$\Delta_{234}$}}

\thicklines 

\put(4,10){\vector(1,0){12}}
\put(24,10){\vector(1,0){12}}
\put(56,10){\vector(-1,0){12}}

\put(46,20){\vector(-1,0){12}}

\put(27,18){\vector(-1,-1){6}}
\put(42,12){\vector(1,1){6}}
\put(31,22){\vector(1,1){6}}

\put(17,12){\vector(-1,1){5}}
\put(22,27){\vector(1,-1){5}}

\qbezier(3,8)(20,2)(37,8)
\put(4.5,7.6){\vector(-3,1){2}}

\end{picture}

\end{center}
\caption{The quivers corresponding to the pseudoline arrangements
  shown in
  Figure~\ref{fig:4-pseudo}. 
The chambers of an arrangement correspond to the 
vertices of the associated quiver. 
}
\label{fig:two-quivers}
\end{figure}

Quiver mutation can be viewed as a generalization of the notion of a
local move: 
there is a combinatorial rule associating a quiver with an
arbitrary pseudoline arrangement so that 
local moves translate into quiver mutations. 
Rather than stating this rule precisely,
we refer to Figure~\ref{fig:two-quivers}, and let the reader guess. 

Let us now define cluster exchanges using the language of quivers. 
This turns out to be very simple. 
Consider a quiver~$Q$ accompanied by an extended cluster~$\zz$,
a~finite collection of algebraically independent elements in our
ambient field of rational functions~$\FFF$.
(Such a pair $(Q,\zz)$ is called a \emph{seed}.) 
Pick a mutable vertex labeled by a cluster variable~$z$.
A~\emph{seed mutation} at~$z$ replaces
$(Q,\zz)$ by the seed $(Q',\zz')$ 
whose quiver is $Q'=\mu_z(Q)$ 
and whose extended cluster is 
$\zz'=\zz\cup\{z'\}\setminus\{z\}$; 
here the new cluster variable~$z'$ is determined by the 
\emph{exchange relation}
\[
z\,z'=
\prod_{z\leftarrow y} y
+\prod_{z\rightarrow y
} y\,.
\]
(The products are over the edges directed at/from~$z$, respectively.) 
For example, the exchange relation associated with the quiver mutation
shown in Figure~\ref{fig:quiver-mutation} is $zz'=vx+uy$;
applying mutation $\mu_x$ to the quiver on the right would invoke the
exchange relation $xx'=z'+u^2$. 

Following the blueprint outlined earlier,
we now define a \emph{cluster algebra} $\AAA(Q)$ 
associated to an arbitrary quiver~$Q$. 
Assign a formal variable to each vertex of~$Q$;
these variables form the initial extended cluster~$\zz$, and 
generate the ambient field~$\FFF$. 
Starting with the initial seed $(Q,\zz)$, repeatedly apply 
seed mutations in all possible directions.
The cluster algebra $\AAA(Q)$ is defined as the subring of $\FFF$
generated by all the elements of all extended clusters obtained by
this recursive process.  

Returning to our running example, we illustrate this definition by
describing the cluster algebra structure in 
$\CC[\operatorname{SL}_4/N]$. 
Let us start with the quiver shown on the left in
Figure~\ref{fig:two-quivers}. 
We view the 9 variables $\Delta_I$ labeling the vertices of this quiver as
formal indeterminates (secretly, they are chamber
minors). 
We declare the variables $\Delta_2$, $\Delta_3$, and~$\Delta_{23}$
mutable; the remaining six variables are frozen. 
There are three possible mutations out of this seed;
we use the quiver to write the corresponding exchange relations:
\begin{align*}
\Delta_2 \Delta_{13}&={\Delta_{12}}\,\Delta_3+{\Delta_1}\,\Delta_{23}\,,\\
\Delta_3\,\Delta_{24}&=
{\Delta_{4}}\,\Delta_{23}+{\Delta_{34}}\,\Delta_{2}\,,\\
\Delta_{23}\,\,\Omega\,\,\,&=
{\Delta_{123}\,\Delta_{34}}\,\Delta_2+{\Delta_{12}\,\Delta_{234}}\,\Delta_3\,. 
\end{align*}
At this point, these relations merely \emph{define} 
$\Delta_{13}$, $\Delta_{24}$, and~$\Omega$ as rational functions in
the original extended cluster. 
The first two relations look familiar: they correspond to
the two local moves that can be applied to the given pseudoline
arrangement. 
The third relation is new: it enables us to flip 
the chamber minor~$\Delta_{23}$, something~we could not do before. 
Although the resulting cluster does not correspond to a pseu\-do\-line
arrangement, 
we can still determine its associated quiver using the definition of
quiver mutation. Continuing this process recursively
\emph{ad infinitum} yields more and more extended clusters;
taken together, they generate a cluster algebra.

If one interprets the elements of the initial cluster as actual flag
minors, 
then the generators produced by this process become rational functions on
the base affine space. 
Remarkably, all these generators are regular functions, 
and generate the ring of all such functions.
This holds for any~$n$, resulting in a cluster algebra structure in
$\CC[\operatorname{SL}_n/N]$; see, e.g.,~\cite{ca3,gls-survey,leclerc-icm}. 

In the special case $n=4$, 
this recursive process produces a \emph{finite} number
of distinct extended clusters, 14~of them to be exact. 
Altogether they contain 15 generators: 
in addition to the $2^4-2=14$ flag minors~$\Delta_I$,  
there is a single new cluster variable 
\[
\Omega=-\Delta_1\Delta_{234}+\Delta_2\Delta_{134}
\]
that already appeared in the third exchange relation above. 

Figure~\ref{fig:sl_4/N} shows the 14 clusters 
for $\CC[\operatorname{SL}_4/N]$ as vertices of a planar graph;  
note that there is one additional vertex at infinity,
so that the graph should be viewed as drawn on a sphere rather than a
plane. 
The regions are labeled by cluster variables. 
Each cluster consists of the three elements labeling the regions
adjacent to the corresponding vertex. 
The edges of the graph correspond to seed mutations. 
The 6 frozen variables are not~shown. 

\begin{figure}[ht]
\begin{center}
\setlength{\unitlength}{0.8pt} 
\begin{picture}(200,200)(0,-6) 
\thinlines 

\put(10,10){\line(-1,-1){10}}
\put(40,80){\line(1,2){20}}
\put(40,80){\line(1,-1){30}}
\put(70,50){\line(3,2){30}}
\put(70,50){\line(1,-4){10}}
\put(80,10){\line(1,0){40}}
\put(60,120){\line(8,-3){40}}
\put(100,70){\line(0,1){35}}
\put(100,70){\line(3,-2){30}}
\put(120,10){\line(1,4){10}}
\put(130,50){\line(1,1){30}}
\put(100,105){\line(8,3){40}}
\put(100,175){\line(0,1){15}}
\put(140,120){\line(1,-2){20}}
\put(190,10){\line(1,-1){10}}

\qbezier(10,10)(10,50)(40,80)
\qbezier(10,10)(40,-5)(80,10)
\qbezier(60,120)(65,155)(100,175)

\qbezier(190,10)(190,50)(160,80)
\qbezier(190,10)(160,-5)(120,10)
\qbezier(140,120)(135,155)(100,175)

\put(10,10){\circle*{4}} 
\put(40,80){\circle*{4}} 
\put(60,120){\circle*{4}} 
\put(70,50){\circle*{4}} 
\put(80,10){\circle*{4}} 
\put(100,70){\circle*{4}} 
\put(100,105){\circle*{4}} 
\put(100,175){\circle*{4}} 
\put(120,10){\circle*{4}} 
\put(130,50){\circle*{4}} 
\put(140,120){\circle*{4}} 
\put(160,80){\circle*{4}} 
\put(190,10){\circle*{4}} 

\put(10,120){\makebox(0,0){$\Delta_{124}$}}
\put(190,120){\makebox(0,0){$\Delta_{134}$}}
\put(100,140){\makebox(0,0){$\Omega$}}
\put(100,40){\makebox(0,0){$\Delta_{23}$}}
\put(100,-17){\makebox(0,0){$\Delta_{14}$}}
\put(73,85){\makebox(0,0){$\Delta_{2}$}}
\put(127,85){\makebox(0,0){$\Delta_{3}$}}
\put(40,30){\makebox(0,0){$\Delta_{24}$}}
\put(160,30){\makebox(0,0){$\Delta_{13}$}}


\end{picture}

\end{center}
\caption{Clusters in $\CC[\operatorname{SL}_4/N]$}
\label{fig:sl_4/N}
\end{figure}

\pagebreak[3]

What do we gain by introducing a cluster algebra structure into a
commutative ring that already appears well understood? 
One reason has been given earlier:
such a structure gives rise to a well-defined notion of the (totally)
positive part of the associated algebraic variety. 
Another reason has to do with defining a ``canonical basis'' in the
algebra at hand; the next paragraph hints at a possible approach. 

Let us call two generators of a cluster algebra \emph{compatible} if
they appear together in some extended cluster. 
A \emph{cluster monomial} is a product of pairwise compatible (not
necessarily distinct) generators. 
It is not too hard to show that in the cluster algebra 
$\AAA=\CC[\operatorname{SL}_4/N]$, 
the cluster monomials form a linear basis.
This is a particular instance of the \emph{dual canonical basis}
of G.~Lusztig (called the ``upper global basis'' by M.~Kashiwara). 

Unfortunately, the general picture (for
arbitrary~$\operatorname{SL}_n$) is much more complicated: 
the cluster monomials seem to form just a part of the dual canonical (or
dual semicanonical) basis; see~\cite{leclerc-icm}. 
The challenge of describing the rest of the dual canonical basis in 
concrete terms remains unmet. 

\medskip

Many other algebraic varieties of 
representation-theoretic importance turn out to possess a natural
structure of a cluster algebra (hence the notions of positivity,
cluster monomials, perhaps canonical bases, etc.).
The list includes Grassmannians, flag manifolds, Schubert varieties,
and double Bruhat cells in arbitrary semisimple Lie groups. 
See \cite{fz-dbc, cdm, gls-survey, keller-bourbaki,leclerc-icm, 
  zelevinsky-millenium, zelevinsky-sf2001}.

\medskip

We conclude this section by mentioning some of the most basic
structural results in the general theory of cluster algebras. 
The first such result is the \emph{Laurent phenomenon}:
the cluster variables are not merely rational functions in the
elements of the initial extended cluster---all of them are in fact Laurent
polynomials! 
We conjectured~\cite{ca1} that these Laurent polynomials always have
positive coefficients; many instances of this conjecture have been
proved (see in particular \cite{carroll-speyer,k-pdf-pos,msw,nakajima})
but the general case seems out of reach at the moment. 
\pagebreak[3]

Another basic structural result is the classification~\cite{ca2} 
of the cluster algebras of \emph{finite type}, i.e., those with finitely many
seeds (equivalently, finitely many generators). 
In the generality presented here, the classification theorem states
that a cluster algebra has finite type if and only if one of its seeds
has a quiver whose subquiver formed by the mutable vertices is an
orientation of a disjoint union of simply-laced Dynkin diagrams. 
(The full-blown  version of the cluster theory leads to a
complete analogue of the Cartan-Killing classification.) 
\pagebreak[3]


The combinatorial scaffolding for a cluster algebra is provided by its
\emph{cluster complex}, a simplicial complex whose
vertices are the cluster variables, and whose maximal simplices are
the clusters. 
In the finite type case, this simplicial complex can be identified as
the dual complex of a \emph{generalized associahedron},
a remarkable convex polytope \cite{cfz, yga} associated with
the corresponding root system. 
In particular, the cluster complex of finite type is homeomorphic to a sphere. 
This can be observed in our running example of
$\CC[\operatorname{SL}_4/N]$: 
the cluster complex is the dual simplicial complex 
of the spherical cell complex shown in Figure~\ref{fig:sl_4/N}. 

\section{Triangulations and laminations
}
\label{sec:teich}

Cluster algebras owe much of their appeal to the ubiquity of the 
combinatorial and algebraic dynamics that underlies them. 
\emph{A~priori}, one might not expect the fairly rigid axioms
governing quiver mutations and exchange relations to be satisfied in a
large variety of contexts. 
Yet this is exactly what happens. 
Moreover, in each instance the framework of clusters and mutations
seems to arise organically rather than artificially. 
A~case in point is discussed in this section: the
classical (by now) machinery of triangulations and laminations on
bordered Riemann surfaces, which goes back to W.~Thurston, 
can be naturally recast in the language of quiver mutations. 
The resulting connection between combinatorial topology and
cluster theory is bound to benefit both. 

This section is based on the papers \cite{cats1, cats2},
which were in turn inspired by the work of 
V.~Fock and A.~Goncharov~\cite{fg-ihes, fock-gonch},
M.~Gekhtman, M.~Shapiro, and A.~Vainshtein~\cite{gsv1, gsv2},
and R.~Penner\cite{penner-lambda}. 

\medskip

\begin{figure}[ht]
\begin{center}
\scalebox{0.5}
{
\epsfbox{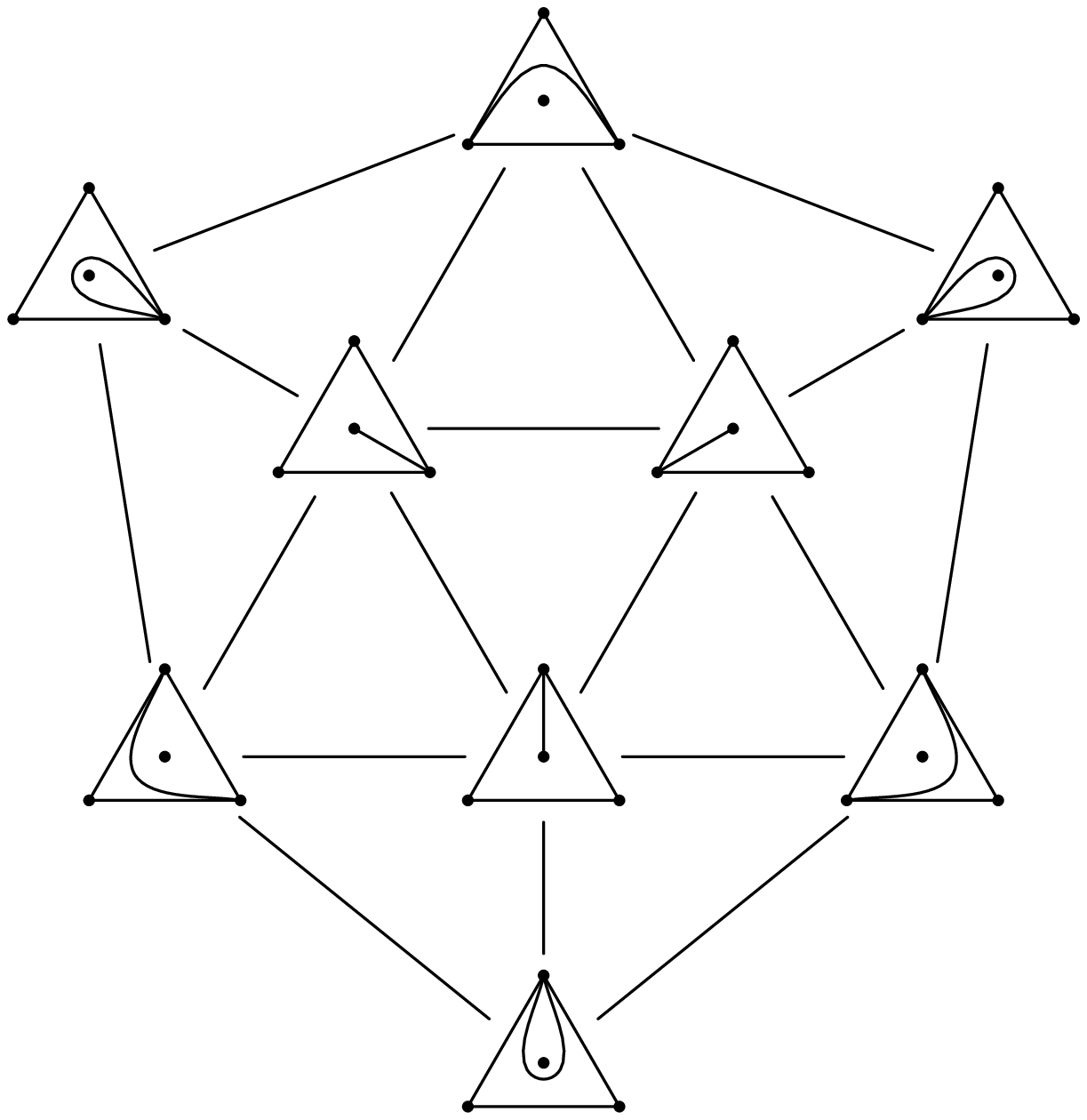}
}
\end{center}
\vspace{-.1in}
\caption{The arc complex of a once-punctured triangle.
Its 10 two-dimensional simplices correspond to ideal triangulations. 
Among them, 6~contain self-folded triangles.}
\label{fig:arc-complex-D3}
\end{figure}

Let $\Surf$ be a connected oriented 
surface with boundary.
(A~few simple 
cases must be ruled out.) 
Fix a finite nonempty set $\Mark$
of \emph{marked points} in the closure of~$\Surf$.
An \emph{arc} in $\SM$ is a non-selfintersecting 
curve in~$\Surf$, considered up to isotopy, 
which connects two points in~$\Mark$, 
does not pass through~$\Mark$, 
and does not cut out an unpunctured monogon or digon. 
Arcs are \emph{compatible} if they have non-intersecting realizations.
Collections of pairwise compatible arcs are the simplices of the 
\emph{arc complex} of~$\Surf$. 
The facets of this simplicial complex correspond to (ideal)
\emph{triangulations}. 
Note that these triangulations may contain \emph{self-folded
  triangles}. 
See Figure~\ref{fig:arc-complex-D3}.

The vertices of the dual graph of the arc complex 
correspond to the triangulations; the edges in this graph correspond
to \emph{flips}. \
A~flip replaces an arc in a triangulation by another (uniquely
defined) arc. 
Note that an edge inside a self-folded triangle cannot be flipped. 
The situation is akin to pseudoline arrangements,
which are likewise related to each other by flips (of a different kind). 

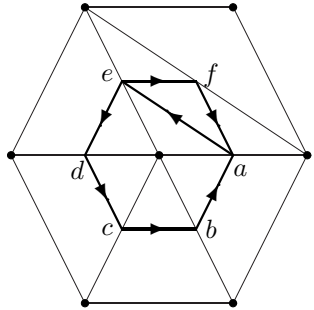
\begin{figure}[ht]
\begin{center}  
\setlength{\unitlength}{2.8pt}  
\begin{picture}(40,40)(0,2)  
\thinlines  
  \put(10,0){\line(-1,2){10}}  
  \put(40,20){\line(-1,2){10}}  
  \put(10,0){\line(1,0){20}}  
  \put(10,40){\line(1,0){20}}  
  \put(0,20){\line(1,2){10}}  
  \put(30,0){\line(1,2){10}}  
  \put(10,0){\line(1,2){10}}  
  \put(0,20){\line(1,0){40}}  
  \put(10,40){\line(1,-2){10}}  
  \put(20,20){\line(1,-2){10}}  
  \put(10,40){\line(3,-2){30}}  

  \put(0,20){\circle*{1}}  
  \put(10,0){\circle*{1}}  
  \put(10,40){\circle*{1}}  
  \put(30,0){\circle*{1}}  
  \put(30,40){\circle*{1}}  
  \put(40,20){\circle*{1}}  
  \put(20,20){\circle*{1}}  

\thicklines  

\put(31,18){\makebox(0,0){$a$}}
\put(27,10){\makebox(0,0){$b$}}
\put(13,10){\makebox(0,0){$c$}}
\put(9,18){\makebox(0,0){$d$}}
\put(13,30.8){\makebox(0,0){$e$}}
\put(27,30.8){\makebox(0,0){$f$}}


\darkred{\put(10,20){\line(1,2){5}}}
\darkred{\put(10,20){\line(1,-2){5}}}
\darkred{\put(15,10){\line(1,0){10}}}
\darkred{\put(25,10){\line(1,2){5}}}
\darkred{\put(30,20){\line(-1,2){5}}}
\darkred{\put(30,20){\line(-3,2){15}}}
\darkred{\put(15,30){\line(1,0){10}}}

\darkred{\put(12,24){\vector(-1,-2){0}}}
\darkred{\put(28,16){\vector(1,2){0}}}
\darkred{\put(13,14){\vector(1,-2){0}}}
\darkred{\put(28,24){\vector(1,-2){0}}}
\darkred{\put(21,10){\vector(1,0){0}}}
\darkred{\put(21,30){\vector(1,0){0}}}
\darkred{\put(21,26){\vector(-3,2){0}}}

\end{picture}  
\quad

\end{center}  
\caption{A triangulation $T$ of a once-punctured hexagon and the
  associated quiver $Q(T)$.}
\label{fig:triang-D6}
\end{figure}

This analogy goes much deeper than it might appear at
first. To see that, we translate the setting into 
the \emph{lingua franca} of quivers. 
Let us define the quiver~$Q(T)$ 
associated to a triangulation~$T$. 
The vertices of $Q(T)$ are labeled by the arcs~in~$T$. \linebreak[3]
If~two arcs belong to the same triangle, 
we connect the corresponding vertices of
the quiver $Q(T)$ by an edge 
whose orientation is determined by the clockwise orientation of the
boundary of the triangle. See Figure~\ref{fig:triang-D6}. 
For triangulations containing self-folded triangles, the definition is
more complicated but is nevertheless completely elementary
and explicit. 


As the reader may have guessed by now,
flips in ideal triangulations translate into mutations of the
associated quivers. 
Furthermore, the quiver language suggests what we should do about the
``forbidden'' flips (of interior edges in self-folded triangles):
forget about triangulations and just mutate the corresponding
quivers.

It is easy to check that a quiver mutation corresponding to 
an edge inside a self-folded triangle transforms any quiver into an
isomorphic one. 
Another simple observation is that the number of different 
(up to isomorphism) quivers $Q(T)$ \linebreak[3]
associated to  
triangulations $T$ of a given surface is \emph{finite}
(because the action of the mapping class group 
on triangulations has finitely many orbits). 
Combining these two observations, one concludes that any quiver~$Q(T)$
associated to a triangulated surface is of \emph{finite mutation
  type}: 
its iterated mutations produce finitely many distinct
(non-isomorphic) quivers. 
In fact, as shown in~\cite{felikson-shapiro-tumarkin}, 
all connected quivers of finite mutation type, with a few exceptions, 
are of the form $Q(T)$, for some triangulation~$T$ of
some marked bordered surface~$\SM$. (We assume that there are no
frozen vertices.) The complete list of exceptions consists of 
(a) quivers with two vertices and more than one edge, and
(b) 11~quivers listed in~\cite{derksen-owen}. 

\begin{figure}[ht]
\begin{center}  
\setlength{\unitlength}{1.8pt}  
\begin{picture}(80,40)(-10,5)  
\thinlines  
\put(0,0){\line(1,0){80}}  
\put(0,40){\line(1,0){80}}  
\put(0,0){\line(0,1){40}}  
\put(80,0){\line(0,1){40}}  

  \put(0,0){\circle*{1}}  
  \put(0,40){\circle*{1}}  
  \put(80,0){\circle*{1}}  
  \put(80,40){\circle*{1}}  

  \put(40,20){\circle{10}}  
\multiput(35,20)(0.5,0){21}{\circle*{0.2}}
\multiput(35.5,20.5)(0.5,0){19}{\circle*{0.2}}
\multiput(35.5,21)(0.5,0){19}{\circle*{0.2}}
\multiput(35.5,21.5)(0.5,0){19}{\circle*{0.2}}
\multiput(36,22)(0.5,0){17}{\circle*{0.2}}
\multiput(36,22.5)(0.5,0){17}{\circle*{0.2}}
\multiput(36.5,23)(0.5,0){15}{\circle*{0.2}}
\multiput(37,23.5)(0.5,0){13}{\circle*{0.2}}
\multiput(37.5,24)(0.5,0){11}{\circle*{0.2}}
\multiput(38,24.5)(0.5,0){9}{\circle*{0.2}}

\multiput(35.5,19.5)(0.5,0){19}{\circle*{0.2}}
\multiput(35.5,19)(0.5,0){19}{\circle*{0.2}}
\multiput(35.5,18.5)(0.5,0){19}{\circle*{0.2}}
\multiput(36,18)(0.5,0){17}{\circle*{0.2}}
\multiput(36,17.5)(0.5,0){17}{\circle*{0.2}}
\multiput(36.5,17)(0.5,0){15}{\circle*{0.2}}
\multiput(37,16.5)(0.5,0){13}{\circle*{0.2}}
\multiput(37.5,16)(0.5,0){11}{\circle*{0.2}}
\multiput(38,15.5)(0.5,0){9}{\circle*{0.2}}

  \put(40,25){\circle*{1}}  

  \put(20,20){\circle*{1}}  
  \put(60,20){\circle*{1}}  

\thicklines  
\darkblue{\qbezier(10,0)(5,20)(10,40)}
\darkblue{\qbezier(12,0)(7,20)(12,40)}
\darkblue{\qbezier(14,0)(-2,55)(80,30)}
\darkblue{\qbezier(13.6,20)(15,45)(80,28)}

\darkblue{\qbezier(40,9)(51,9)(55,20)}
\darkblue{\qbezier(40,9)(29,9)(25,20)}
\darkblue{\qbezier(55,20)(56.5,24)(60,24)}
\darkblue{\qbezier(25,20)(23.5,24)(20,24)}
\darkblue{\qbezier(63,20)(63,24)(60,24)}
\darkblue{\qbezier(17,20)(17,24)(20,24)}
\darkblue{\qbezier(63,20)(63,17)(60,17)}
\darkblue{\qbezier(17,20)(17,17)(20,17)}
\darkblue{\qbezier(13.6,20)(13.6,14)(18,14)}
\darkblue{\qbezier(20,17)(22,17)(22,20)}
\darkblue{\qbezier(60,17)(58,17)(58,20)}
\darkblue{\qbezier(18,14)(23,14)(23,20)}
\darkblue{\qbezier(23,20)(23,22.5)(20,22.5)}
\darkblue{\qbezier(20,21.5)(22,21.5)(22,20)}
\darkblue{\qbezier(60,21.5)(58,21.5)(58,20)}
\darkblue{\qbezier(20,21.5)(18.7,21.5)(18.7,20)}
\darkblue{\qbezier(60,21.5)(61.3,21.5)(61.3,20)}
\darkblue{\qbezier(20,22.5)(19,22.5)(18.6,22)}

  \darkblue{\put(40,20){\circle{14}}}  
  \darkblue{\put(40,20){\circle{18}}}

\end{picture}  
\setlength{\unitlength}{3.6pt}
\begin{picture}(40,20)(-5,2.5)
\thinlines

  \put(20,20){\circle*{1}}
  \put(15,10){\circle*{1}}

  \put(20,10){\circle{20}}

\thicklines
\lightred{\qbezier(15,18.4)(20,15)(25,18.4)}
\lightred{\qbezier(15,1.6)(20,5)(25,1.6)}

\lightred{\put(15,10){\circle{6}}}
\lightred{\put(25,10){\circle{6}}}

\end{picture}
\end{center}
\caption{(a) A lamination; (b) curves not allowed in a lamination.}
\label{fig:lamin-yes-no}
\end{figure}
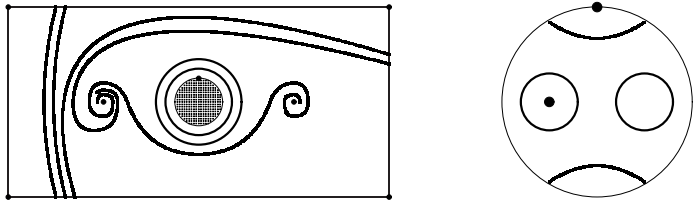

The construction of quivers $Q(T)$
can be generalized by involving
W.~Thurston's machinery of laminations on Riemann surfaces. 
An integral (unbounded measured) \emph{lamination} on~$\SM$ 
is a finite collection of non-selfintersecting and
pairwise non-inter\-secting curves in~$\Surf$,  
considered modulo isotopy. 
The curves in a lamination must satisfy certain constraints.
In particular, 
each of them is either closed, or runs from boundary to boundary, or
spirals into an interior marked point (a puncture). 
See Figure~\ref{fig:lamin-yes-no}. 

\begin{figure}[ht]
\begin{center}
\setlength{\unitlength}{1pt}
\begin{picture}(60,60)(0,2)
\thicklines
  \put(0,20){\line(1,2){20}}
  \put(0,20){\line(1,-1){20}}
  \put(20,0){\line(0,1){60}}
  \put(20,0){\line(1,1){40}}
  \put(20,60){\line(2,-1){40}}

  \put(20,0){\circle*{3}}
  \put(20,60){\circle*{3}}
  \put(0,20){\circle*{3}}
  \put(60,40){\circle*{3}}

\put(24,43){\makebox(0,0){$\gamma$}}
\put(-15,35){\makebox(0,0){$+1$}}

\darkblue{\qbezier(-10,10)(30,20)(60,60)}

\end{picture}
\qquad\qquad\qquad
\begin{picture}(60,60)(0,2)
\thicklines
  \put(0,20){\line(1,2){20}}
  \put(0,20){\line(1,-1){20}}
  \put(20,0){\line(0,1){60}}
  \put(20,0){\line(1,1){40}}
  \put(20,60){\line(2,-1){40}}

  \put(20,0){\circle*{3}}
  \put(20,60){\circle*{3}}
  \put(0,20){\circle*{3}}
  \put(60,40){\circle*{3}}

\put(-15,35){\makebox(0,0){$-1$}}
\put(16,20){\makebox(0,0){$\gamma$}}

\darkblue{\qbezier(-10,60)(30,40)(60,10)}

\end{picture}
\end{center}
\caption{A (signed) contribution of a curve in $L$ to the shear
  coordinate~$b_\gamma(T,L)$.} 
\label{fig:shear-sign}
\end{figure}
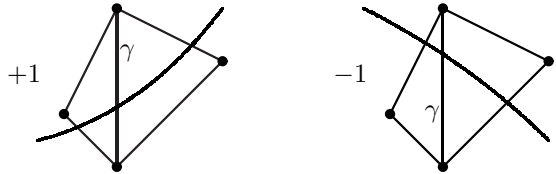

Let $L$ be an integral lamination, and $T$ a triangulation 
without self-folded triangles. 
For an arc $\gamma$ in~$T$, the \emph{shear coordinate} $b_\gamma(T,L)$ is
the signed number of curves in~$L$ which
intersect~$\gamma$ and in doing so, connect the opposite sides of the
quadrilateral surrounding~$\gamma$.
The sign depends on which pair of opposite sides the curves 
connect;
see Figure~\ref{fig:shear-sign}. 

\pagebreak[3]

By a theorem of W.~Thurston, the shear coordinates \emph{coordinatize}
integral laminations in the following sense: 
for a fixed triangulation~$T$, 
the map 
\[
L\mapsto (b_\gamma(T,L))_{\gamma\in T}
\]
is a bijection between 
integral 
laminations and~$\ZZ^n$. 

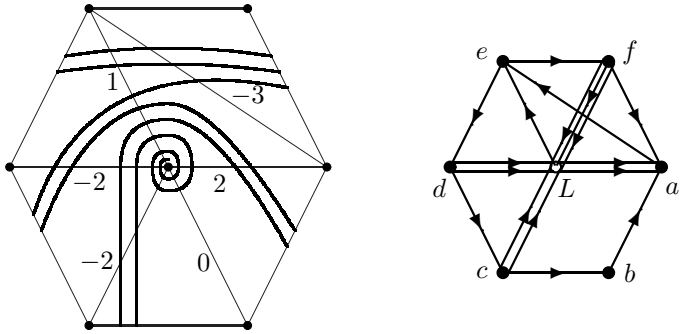
\begin{figure}[ht]
\vspace{-.4in}
\begin{center}  
\setlength{\unitlength}{3pt}  
\begin{picture}(40,41)(0,0)  
\thinlines  
  \put(10,0){\line(-1,2){10}}  
  \put(40,20){\line(-1,2){10}}  
  \put(10,0){\line(1,0){20}}  
  \put(10,40){\line(1,0){20}}  
  \put(0,20){\line(1,2){10}}  
  \put(30,0){\line(1,2){10}}  
  \put(10,0){\line(1,2){10}}  
  \put(0,20){\line(1,0){40}}  
  \put(10,40){\line(1,-2){10}}  
  \put(20,20){\line(1,-2){10}}  
  \put(10,40){\line(3,-2){30}}  

  \put(0,20){\circle*{1}}  
  \put(10,0){\circle*{1}}  
  \put(10,40){\circle*{1}}  
  \put(30,0){\circle*{1}}  
  \put(30,40){\circle*{1}}  
  \put(40,20){\circle*{1}}  
  \put(20,20){\circle*{1}}  

\thicklines  

\darkblue{\qbezier(16,0)(16,10)(16,20)}
\darkblue{\qbezier(16,20)(16,24)(20,24)}
\darkblue{\qbezier(23,20)(23,24)(20,24)}
\darkblue{\qbezier(23,20)(23,17)(20,17)}
\darkblue{\qbezier(20,17)(18,17)(18,20)}
\darkblue{\qbezier(20,22)(18,22)(18,20)}
\darkblue{\qbezier(20,22)(21.3,22)(21.3,20)}
\darkblue{\qbezier(20,18.5)(21.3,18.5)(21.3,20)}
\darkblue{\qbezier(20,18.5)(19,18.5)(19,20)}
\darkblue{\qbezier(20,21)(19,21)(19,20)}
 
\darkblue{\qbezier(14,0)(14,10)(14,20)}
\darkblue{\qbezier(14,20)(14,26)(20,26)}

\darkblue{\qbezier(35,10)(26,26)(20,26)}

\darkblue{\qbezier(36,12)(26,28)(20,28)}
\darkblue{\qbezier(4,12)(10,28)(20,28)}
\darkblue{\qbezier(3,14)(10,35)(35,30)}

\darkblue{\qbezier(6,32)(20,34)(34,32)}
\darkblue{\qbezier(7,34)(20,36)(33,34)}

\darkblue{
\put(24.5,8){\makebox(0,0){$0$}}
\put(11,8){\makebox(0,0){$-2$}}
\put(10,18){\makebox(0,0){$-2$}}
\put(26.5,18){\makebox(0,0){$2$}}
\put(30,29){\makebox(0,0){$-3$}}
\put(13,30.8){\makebox(0,0){$1$}}
}
\end{picture}  
\setlength{\unitlength}{4pt}  
\begin{picture}(40,40)(0,5)  
\thicklines  

\put(31,18){\makebox(0,0){$a$}}
\put(27,10){\makebox(0,0){$b$}}
\put(13,10){\makebox(0,0){$c$}}
\put(9,18){\makebox(0,0){$d$}}
\put(13,30.8){\makebox(0,0){$e$}}
\put(27,30.8){\makebox(0,0){$f$}}

\put(21,18){\makebox(0,0){$L$}}

\darkred{\put(10,20){\circle*{1}}}
\darkred{\put(15,10){\circle*{1}}}
\darkred{\put(15,30){\circle*{1}}}
\darkred{\put(25,10){\circle*{1}}}
\darkred{\put(30,20){\circle*{1}}}
\darkred{\put(25,30){\circle*{1}}}

\darkred{\put(20,20){\circle{1}}}

\darkred{\put(10,20){\line(1,2){5}}}
\darkred{\put(10,20){\line(1,-2){5}}}
\darkred{\put(15,10){\line(1,0){10}}}
\darkred{\put(25,10){\line(1,2){5}}}
\darkred{\put(30,20){\line(-1,2){5}}}
\darkred{\put(30,20){\line(-3,2){15}}}
\darkred{\put(15,30){\line(1,0){10}}}

\darkred{\put(12,24){\vector(-1,-2){0}}}
\darkred{\put(28,16){\vector(1,2){0}}}
\darkred{\put(13,14){\vector(1,-2){0}}}
\darkred{\put(28,24){\vector(1,-2){0}}}
\darkred{\put(21,10){\vector(1,0){0}}}
\darkred{\put(21,30){\vector(1,0){0}}}
\darkred{\put(18,28){\vector(-3,2){0}}}

\darkred{\put(10.2,20.4){\line(1,0){9.6}}}
\darkred{\put(10.2,19.6){\line(1,0){9.6}}}
\darkred{\put(17,20.4){\vector(1,0){0}}}
\darkred{\put(17,19.6){\vector(1,0){0}}}

\darkred{\put(20.2,20.4){\line(1,0){9.6}}}
\darkred{\put(20.2,19.6){\line(1,0){9.6}}}
\darkred{\put(27,20.4){\vector(1,0){0}}}
\darkred{\put(27,19.6){\vector(1,0){0}}}

\darkred{\put(14.6,10.2){\line(1,2){5}}}
\darkred{\put(17.6,16.2){\vector(1,2){0}}}
\darkred{\put(15.4,9.8){\line(1,2){5}}}
\darkred{\put(18.4,15.8){\vector(1,2){0}}}

\darkred{\put(25.5,29.75){\line(-1,-2){5}}}
\darkred{\put(21.5,21.75){\vector(-1,-2){0}}}
\darkred{\put(24.5,30.25){\line(-1,-2){5}}}
\darkred{\put(20.5,22.25){\vector(-1,-2){0}}}
\darkred{\put(24.9,29.8){\line(-1,-2){4.8}}}
\darkred{\put(22.9,25.8){\vector(-1,-2){0}}}

\darkred{\put(15.1,29.8){\line(1,-2){4.8}}}
\darkred{\put(17,26){\vector(-1,2){0}}}

\end{picture}  
\end{center}  
\vspace{-.1in}
\caption{(a) Shear coordinates of a lamination~$L$; 
(b) the quiver $Q(T,\{L\})$.}
\label{fig:lamin-shear}
\end{figure}

A \emph{multi-lamination} $\LL$ on~$\SM$ is an arbitrary finite family of
laminations. 
Given such $\LL$ and a triangulation~$T$ of the surface~$\SM$,
we construct the ``extended'' quiver $Q(T,\LL)$ by adding vertices and
oriented edges to~$Q(T)$ as follows. 
For each lamination $L$ in~$\LL$, we introduce a new vertex
labeled~by~$L$. We then connect this vertex to each vertex in~$Q$, 
say labeled by an arc~$\gamma$, by $|b_\gamma(T,L)|$ 
edges whose direction is determined by the sign of~$b_\gamma(T,L)$. 
See Figure~\ref{fig:lamin-shear}. 

Amazingly, the same property as before holds: 
for a fixed multi-lamination~$\LL$, 
a flip in a triangulation~$T$ translates into the corresponding
mutation in the quiver $Q(T,\LL)$.
(The definition of the latter can be generalized to allow for self-folded
triangles.)
This strongly suggests the existence of a cluster algebra structure
associated with any given marked surface~$\SM$ and any 
multi-lamination~$\LL$~on~it. 

This class of cluster algebras can be 
understood on several levels. 
On the combinatorial level, the
cluster complex of such an algebra can be explicitly described 
in terms of \emph{tagged arcs}, which are ordinary arcs 
adorned with very simple combinatorial decorations. 
This description represents the cluster complex
as a finite covering space for the arc complex. 
The 
cluster complex turns out to be either contractible or 
homotopy equivalent to a sphere. 
Unlike the generalized associahedra mentioned above,
these cluster complexes are usually not compact;
moreover, with a few exceptions, they exhibit exponential growth. 
See~\cite{cats1}. 


The coordinatization theorem implies that any quiver~$Q$ 
whose mutable part can be interpreted as a quiver $Q(T)$ corresponding
to a triangulation~$T$ of some marked surface~$\SM$, 
there exists a (unique) multi-lamination~$\LL$ on~$\SM$ such that
$Q=Q(T,\LL)$. 
In view of the discussion above, 
the cluster algebra $\AAA(Q)$ associated with such a quiver~$Q$
depends only on~$\SM$ and~$\LL$ but not on the triangulation~$T$. 
Consequently, one should be able to understand this cluster algebra
in terms of the topology of the surface $\SM$ and the
multi-lamination~$\LL$.

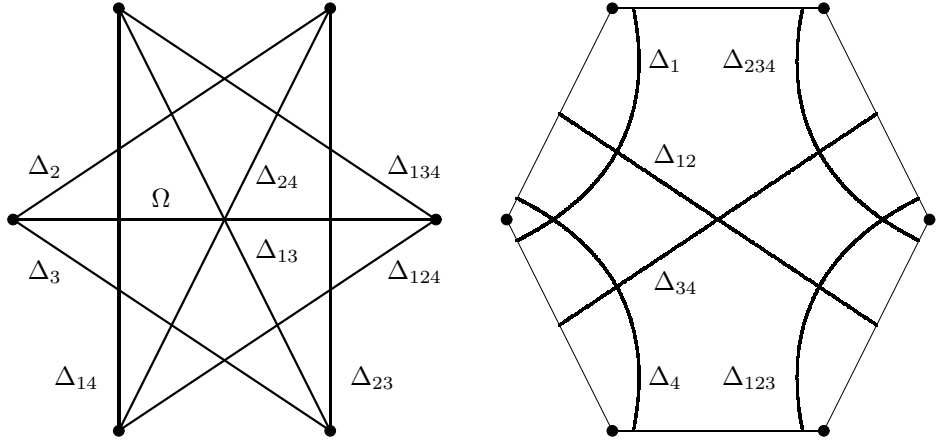
\begin{figure}[ht]
\begin{center}  
\setlength{\unitlength}{4pt}  
\begin{picture}(40,42)(0,0)  
\thinlines  

  \put(0,20){\circle*{1}}  
  \put(10,0){\circle*{1}}  
  \put(10,40){\circle*{1}}  
  \put(30,0){\circle*{1}}  
  \put(30,40){\circle*{1}}  
  \put(40,20){\circle*{1}}  

\thicklines
  \put(0,20){\line(3,2){30}} 
\put(3,25){\makebox(0,0){$\Delta_2$}}
  \put(0,20){\line(3,-2){30}} 
\put(3,15){\makebox(0,0){$\Delta_3$}}
  \put(30,0){\line(0,1){40}} 
\put(34,5){\makebox(0,0){$\Delta_{23}$}}

  \put(40,20){\line(-3,2){30}} 
\put(38,25){\makebox(0,0){$\Delta_{134}$}}
  \put(40,20){\line(-3,-2){30}} 
\put(38,15){\makebox(0,0){$\Delta_{124}$}}
  \put(10,0){\line(0,1){40}} 
\put(6,5){\makebox(0,0){$\Delta_{14}$}}

  \put(10,0){\line(1,2){20}} 
\put(25,24){\makebox(0,0){$\Delta_{24}$}}
  \put(10,40){\line(1,-2){20}} 
\put(25,17){\makebox(0,0){$\Delta_{13}$}}
  \put(0,20){\line(1,0){40}} 
\put(14,22){\makebox(0,0){$\Omega$}}

\end{picture}  
\qquad
\setlength{\unitlength}{4pt}  
\begin{picture}(40,42)(0,0)  
\thinlines  
  \put(10,0){\line(-1,2){10}}  
  \put(40,20){\line(-1,2){10}}  
  \put(10,0){\line(1,0){20}}  
  \put(10,40){\line(1,0){20}}  
  \put(0,20){\line(1,2){10}}  
  \put(30,0){\line(1,2){10}}  

  \put(0,20){\circle*{1}}  
  \put(10,0){\circle*{1}}  
  \put(10,40){\circle*{1}}  
  \put(30,0){\circle*{1}}  
  \put(30,40){\circle*{1}}  
  \put(40,20){\circle*{1}}  

\thicklines  
\darkred{\qbezier(12,0)(15,15)(1,22)
\put(15,5){\makebox(0,0){$\Delta_4$}}
}
\lightblue{\qbezier(12,40)(15,25)(1,18)
\put(15,35){\makebox(0,0){$\Delta_1$}}
}
\darkblue{\qbezier(28,0)(25,15)(39,22)
\put(23,5){\makebox(0,0){$\Delta_{123}$}}
}
\lightred{\qbezier(28,40)(25,25)(39,18)
\put(23,35){\makebox(0,0){$\Delta_{234}$}}
}
\darkgreen{\qbezier(5,10)(20,20)(35,30)
\put(16,14){\makebox(0,0){$\Delta_{34}$}}
}
\lightgreen{\qbezier(5,30)(20,20)(35,10)
\put(16,26){\makebox(0,0){$\Delta_{12}$}}
}

\end{picture}  
\end{center}  
\caption{(a) Labeling the cluster
  variables in $\mathbb{C}[\operatorname{SL}_4/N]$
by the diagonals of a hexagon.
(b)~Labeling the frozen variables 
by laminations, each consisting of a single curve.}
\label{fig:lamin-SL3/N}
\end{figure}

We illustrate this construction by returning, once again, to the
example of the cluster algebra $\AAA=\mathbb{C}[\operatorname{SL}_4/N]$. 
The mutable part of any quiver~$Q$ defining this algebra 
(see, e.g., Figure~\ref{fig:two-quivers}) 
has 3~vertices, and is
isomorphic to a quiver $Q(T)$ associated to a triangulation of a
hexagon, i.e., a disk with 6~marked points on the boundary. 
Thus, we can let $\SM$ be a hexagon. 
Due to the absence of marked points in the interior of~$\Surf$,
the construction simplifies considerably: 
there are no self-folded triangles, and 
the cluster complex coincides with the arc complex. 
The underlying combinatorics of $\AAA$ is thus modeled 
as follows: 
cluster variables correspond to arcs (that is, the diagonals of the hexagon), 
clusters correspond to triangulations, 
and exchanges correspond to flips. 
It remains to determine the appropriate multi-lamination~$\LL$.
This is done by interpreting the multiplicities of edges connecting
the frozen vertices in~$Q$ to the mutable ones as shear coordinates of
laminations, and then constructing the unique laminations having those
shear coordinates. 
The result is shown in Figure~\ref{fig:lamin-SL3/N}. 

It is natural to ask whether cluster variables in the cluster algebra 
associated with a multi-lamination on a bordered surface 
can be given an intrinsic geometric interpretation.  
The answer is yes:
each cluster variable can be viewed as a suitably renormalized 
\emph{lambda length}~\cite{penner-lambda} (a.k.a.\ Penner coordinate) 
of the corresponding (tagged) arc. 
For a given arc, such a lambda length is a real function on 
(an appropriate generalization of) the  
\emph{decorated Teichm\"uller space} for~$\SM$;
see \cite{cats2} for further details. 
Thus in this geometric realization,
the decorated Teichm\"uller space plays the role of the
corresponding totally positive variety. 
This brings us back full circle to the problems discussed at the
end of Section~\ref{sec:tp}, namely to  
the challenges of understanding the stratification of a totally
nonnegative variety
(in this case, a compactified decorated Teichm\"uller space)
and the singularities of its boundary.


\newpage

\begin{thebibliography}{xxx}
\bibitem{ando}
T.~Ando, Totally positive matrices, 
{\sl Linear Algebra Appl.}\ 
{\bf 90} (1987), 165--219.

\bibitem{BFZ} A.~Berenstein, S.~Fomin, and A.~Zelevinsky, 
Parametrizations of canonical bases
and totally positive matrices,
{\sl Adv.\ Math.}\ {\bf 122} (1996), 49--149.

\bibitem{ca3}
A.~Berenstein, S.~Fomin and A.~Zelevinsky,
Cluster algebras~III: Upper bounds and double Bruhat cells, 
\textsl{Duke Math.~J.}\ \textbf{126} (2005), 1--52. 


\bibitem{carroll-speyer}
G.~Carroll and D.~Speyer, 
The cube recurrence, 
\textsl{Electron.\ J.\ Combin.}\ \textbf{11} (2004), no.~1, 
Research Paper~73, 31~pp. 

\bibitem{carter}
R.~W.~Carter, \textsl{Cluster algebras}. 
Textos de Matem\'atica, S\'erie~B, 37. 
Universidade de Coimbra, 2006. 

\bibitem{cfz}
F.~Chapoton, S.~Fomin, and A.~Zelevinsky,
Polytopal realizations of generalized associahedra,
\textsl{Canad.\ Math.\ Bull.}\ \textbf{45} (2002), 537--566. 

\bibitem{chekhov}
L.~O.~Chekhov, 
Orbifold Riemann surfaces and geodesic algebras, 
\textsl{J.~Phys.~A} \textbf{42} (2009), no. 30, 304007, 32 pp. 

\bibitem{cryer}
C.~Cryer,
The $LU$-factorization of totally positive matrices,
{\sl Linear Algebra Appl.}\ \textbf{7} (1973), 83--92.

\bibitem{cryer76}
C.~Cryer,
Some properties of totally positive matrices,
{\sl Linear Algebra Appl.}\ \textbf{15} (1976), 1--25.

\bibitem{derksen-owen}
H.~Derksen and T.~Owen, 
New graphs of finite mutation type, 
\textsl{Electron.\ J.\ Combin.}\ \textbf{15} (2008), no. 1, 
Research Paper 139, 15 pp. 

\bibitem{dfz1}
H.~Derksen, J.~Weyman, and A.~Zelevinsky, 
Quivers with potentials and their representations~I: Mutations, 
\textsl{Selecta Math.\ (N.S.)}  \textbf{14} (2008), 59-119.

\bibitem{dfz2}
H.~Derksen, J.~Weyman, and A.~Zelevinsky, 
Quivers with potentials and their representations~II:
Applications to cluster algebras,
to appear in \textsl{J.~Amer.\ Math.\ Soc.}, \texttt{arXiv:0904.0676}. 

\bibitem{k-pdf}
P.~Di Francesco and R.~Kedem, 
$Q$-systems as cluster algebras.~II. 
Cartan matrix of finite type and the polynomial property, 
\textsl{Lett.\ Math.\ Phys.}\ \textbf{89} (2009), 183--216. 

\bibitem{k-pdf-pos}
P.~Di Francesco and R.~Kedem, 
$Q$-systems, heaps, paths and cluster positivity,
\textsl{Comm.\ Math.\ Phys.}\ \textbf{293} (2010), 727--802.

\bibitem{felikson-shapiro-tumarkin}
A.~Felikson, M.~Shapiro, and P.~Tumarkin,
Skew-symmetric cluster algebras of finite mutation type,
\texttt{arXiv:0811.1703}.

\bibitem{fg-ihes}
V.~V.~Fock and A.~B.~Goncharov, 
Moduli spaces of local systems and higher {T}eichm\"uller
  theory, 
\textsl{Publ.\ Math.\ Inst.\ Hautes \'Etudes Sci.}\ (2006), 
no.~103, 1--211. 

\bibitem{fock-gonch}
V.~V.~Fock and A.~B.~Goncharov, 
Dual Teichmuller and lamination spaces,
\textsl{Hand\-book of Teichm\"uller theory, vol.~I,}  647--684,
Eur.\ Math.\ Soc., Z\"urich, 2007. 

\bibitem{portal}
S.~Fomin, 
Cluster Algebras Portal, \\
\texttt{http://www.math.lsa.umich.edu/$\tilde{\
  }$fomin/cluster.html}. 

\bibitem{pcmi}
S.~Fomin and N.~Reading, 
Root systems and generalized associahedra,
\textsl{Geometric Combinatorics (Park City, UT, 2003),} 63--131, 
IAS/Park City Math.\ Ser., 14, Amer.\ Math.\ Soc., Providence, RI, 2007. 

\bibitem{cats1}
S.~Fomin, M.~Shapiro, and D.~Thurston, 
Cluster algebras and triangulated surfaces. Part I: Cluster complexes,
\textsl{Acta Math.}\ \textbf{201} (2008), 83--146. 

\bibitem{cats2}
S.~Fomin and D.~Thurston, 
Cluster algebras and triangulated surfaces. Part II: Lambda lengths, 
preprint. 

\bibitem{fz-dbc}
S.~Fomin and A.~Zelevinsky,
Double Bruhat cells and total positivity, 
\textsl{J.~Amer.\ Math.\ Soc.}\ \textbf{12} (1999), 335--380. 

\bibitem{fs}
S.~Fomin and M.~Z.~Shapiro, 
Stratified spaces formed by totally positive varieties, 
\textsl{Michigan Math.~J.}\ \textbf{48} (2000), 253--270. 

\bibitem{fz-osc}
S.~Fomin and A.~Zelevinsky,
Totally nonnegative and oscillatory elements in semi\-simple groups,
\textsl{Proc.\ Amer.\ Math.\ Soc.}\ \textbf{128} (2000), 3749--3759. 

\bibitem{tptp}
S.~Fomin and A.~Zelevinsky,
Total positivity: tests and parametrizations,
\textsl{Math.\ Intelligencer} \textbf{22} (2000), 23--33. 

\bibitem{ca1}
S.~Fomin and A.~Zelevinsky,
Cluster algebras I: Foundations, \textsl{J.~Amer.\ Math.\ Soc.}\
\textbf{15} (2002), 497--529. 

\bibitem{ca2}
S.~Fomin and A.~Zelevinsky,
Cluster algebras~II: Finite type classification, 
\textsl{Invent.\ Math.}\ \textbf{154} (2003), 63--121. 

\bibitem{yga}
S.~Fomin and A.~Zelevinsky,
$Y$-systems and generalized associahedra, 
\textsl{Ann.\ of Math.}\ \textbf{158} (2003), 977--1018. 

\bibitem{cdm}
S.~Fomin and A.~Zelevinsky,
Cluster algebras: Notes for the CDM-03 conference,
\textsl{Current Developments in Mathematics, 2003}, 1--34, 
Int.\ Press, 2004. 

\bibitem{ca4}
S.~Fomin and A.~Zelevinsky,
Cluster algebras~IV: Coefficients, \textsl{Compos.\ Math.}\ 
\textbf{143} (2007), 112--164. 

\bibitem{GK}
F.~R.~Gantmacher and M.~G.~Krein, 
{\sl Oscillation matrices and kernels and small vibrations of
  mechanical systems},  
AMS Chelsea Publishing, Providence, RI, 2002.
(Original Russian edition, 1941.) 

\bibitem{gsv1}
M.~Gekhtman, M.~Shapiro, and A.~Vainshtein, 
Cluster algebras and Poisson geometry, 
\textsl{Mosc.\ Math.~J.}\ \textbf{3} (2003), 899--934.

\bibitem{gsv2}
M.~Gekhtman, M.~Shapiro, and A.~Vainshtein, 
Cluster algebras and Weil-{P}etersson forms, 
\textsl{Duke Math.~J.}\ \textbf{127} (2005), 291--311.

\bibitem{gls-survey}
C.~Geiss, B.~Leclerc, and J.~Schr\"oer, 
Preprojective algebras and cluster algebras,
\emph{in:}
Trends in representation theory of algebras and related topics, 253--283,
\textsl{EMS Ser.\ Congr.\ Rep.}, Eur. Math. Soc., Zürich, 2008. 

\bibitem{gasca-micchelli}
{\sl Total positivity and its applications}, M.~Gasca and
C.~A.~Micchelli (Eds.), {\sl Mathematics and its Applications} {\bf
359}, Kluwer Academic Publishers, Dordrecht, 1996.

\bibitem{hersh}
P.~Hersh,
Regular cell complexes in total positivity,
\texttt{arXiv:0711.1348}.

\bibitem{iikns}
R.~Inoue, O.~Iyama, A.~Kuniba, T.~Nakanishi, and J.~Suzuki,
Periodicities of $T$-systems and $Y$-systems,
to appear in \textsl{Nagoya Math.~J.}, 
\texttt{arXiv:0812.0667}. 

\bibitem{karlin} 
S.~Karlin, {\sl Total positivity}, 
Stanford University Press, 1968. 

\bibitem{keller-bourbaki}
B.~Keller,
Alg\`ebres amass\'ees et applications,
\textsl{S\'eminaire Bourbaki, 2009/2010, expos\'e~1014},
\texttt{arXiv:0911.2903}. 

\bibitem{keller-categorification}
B.~Keller,
Categorification of acyclic cluster algebras: an introduction,
to appear in: \textsl{Proceedings of the conference ``Higher structures in
Geometry and Physics 2007,''} Birkh\"auser.

\bibitem{keller-course}
B.~Keller,
Cluster algebras, quiver representations and triangulated categories,
to appear in: \textsl{Triangulated categories,}
London Math.\ Soc., 2010. 


\bibitem{kontsevich-soibelman}
M.~Kontsevich and Y.~Soibelman,
Stability structures, motivic Donaldson-Thomas invariants and cluster
transformations, 
\texttt{arXiv:0811.2435}, November 2008. 

\bibitem{leclerc-icm}
B.~Leclerc,
Cluster algebras and representation theory,
\textsl{Proceedings of the International Congress of Mathematicians}, 
Hyderabad, India, 2010. 

\bibitem{littelmann}
P.~Littelmann,
Bases canoniques et applications,
\textsl{S\'eminaire Bourbaki, 1997/1998, expos\'e~847};
\textsl{Ast\'erisque} \textbf{252} (1998), 287--306. 

\bibitem{loewner}
C.~Loewner, On totally positive matrices,
{\sl Math.~Z.}\ {\bf 63} (1955), 338--340.

\bibitem{lusztig}
G.~Lusztig, Total positivity in reductive groups,
in: {\sl Lie theory and geometry: 
in honor of Bertram Kostant, Progress in Mathematics} {\bf 123}, 
Birkh\"auser, 1994, 531--568.

\bibitem{lusztig-survey}
G.~Lusztig, Introduction to total positivity,
in: {\sl Positivity in Lie theory: open problems,
de Gruyter Exp.\ Math.}\ {\bf 26}, 
de Gruyter, Berlin, 1998, 133--145.

\bibitem{msw}
G.~Musiker, R.~Schiffler, and L.~Williams,
Positivity for cluster algebras from surfaces,
\texttt{arXiv:0906.0748}.

\bibitem{nagao}
K.~Nagao,
Donaldson-Thomas theory and cluster algebras,
\texttt{arXiv:1002.4884}, February 2010. 

\bibitem{nakajima}
H.~Nakajima,
Quiver varieties and cluster algebras,
\texttt{arXiv:0905.0002}, May 2009. 

\bibitem{penner-lambda}
R.~C.~Penner,
Lambda lengths, 
 University of Aarhus, lecture notes, August 2006,
\texttt{http://www.ctqm.au.dk/research/MCS/lambdalengths.pdf}.

\bibitem{reading-speyer}
N.~Reading and D.~Speyer, 
Cambrian fans, 
\textsl{J.~Eur.\ Math.\ Soc.}\ \textbf{11} (2009), 407--447. 

\bibitem{schoenberg}
I.~J.~Schoenberg, 
\"Uber Variationsverminderende lineare Transformationen,~\textsl{Math.~Z.} 
\textbf{32} (1930), 321--328. 

\bibitem{whitney}
A.~M.~Whitney,
A reduction theorem for totally positive matrices,
{\sl J.~d'Analyse Math.}\ {\bf 2} (1952), 88--92.

\bibitem{zelevinsky-millenium}
A.~Zelevinsky, 
Cluster algebras: origins, results and conjectures. 
\textsl{Advances in algebra towards millennium problems}, 
85--105, SAS Int.\ Publ., Delhi, 2005. 

\bibitem{zelevinsky-sf2001}
A.~Zelevinsky, 
From Littlewood-Richardson coefficients to cluster algebras in three
lectures, in: \textsl{Symmetric functions 2001: surveys of
developments and perspectives}, edited by S.~Fomin,  253--273, 
NATO Sci. Ser. II Math. Phys. Chem., 74, Kluwer Acad. Publ., Dordrecht, 2002. 

\bibitem{zelevinsky-whatis}
A.~Zelevinsky, 
What is $\dots$ a cluster algebra?  
\textsl{Notices Amer.\ Math.\ Soc.}\  \textbf{54}  (2007),  no. 11,
1494--1495. 

\end{thebibliography}
\end{document}